\documentclass[final,times]{elsarticle}
\usepackage[top=2.5cm,bottom=2.5cm,left=1.35cm,right=1.35cm]{geometry}

\usepackage{amsmath,amssymb,amsfonts,bm}
\usepackage{graphicx}
\usepackage{booktabs}
\usepackage{subcaption}
\usepackage{placeins}
\usepackage{pgfplots}
\usepackage[hidelinks]{hyperref}
\usepackage{qcircuit}

\usepackage{comment}
\pgfplotsset{compat=1.18}

\graphicspath{{figures}}

\usepackage{soul,color}
\soulregister\cite7
\soulregister\ref7
\soulregister\eqref7
\soulregister\pageref7
\newcommand{\hll}[1]{%
  \colorbox{yellow}{$\displaystyle#1$}}

\sethlcolor{white} 
\renewcommand{\hll}[1]{%
  \colorbox{white}{$\displaystyle#1$}}

\usepackage{lineno}

\newcommand*\patchAmsMathEnvironmentForLineno[1]{%
  \expandafter\let\csname old#1\expandafter\endcsname\csname #1\endcsname
  \expandafter\let\csname oldend#1\expandafter\endcsname\csname end#1\endcsname
  \renewenvironment{#1}%
    {\linenomath\csname old#1\endcsname}%
    {\csname oldend#1\endcsname\endlinenomath}%
}

\newcommand*\patchBothAmsMathEnvironmentsForLineno[1]{%
  \patchAmsMathEnvironmentForLineno{#1}%
  \patchAmsMathEnvironmentForLineno{#1*}%
}

\AtBeginDocument{%
  \patchBothAmsMathEnvironmentsForLineno{equation}
  \patchBothAmsMathEnvironmentsForLineno{align}
  \patchBothAmsMathEnvironmentsForLineno{gather}
  \patchBothAmsMathEnvironmentsForLineno{multline}
}

\usepackage{etoolbox}
\makeatletter
\patchcmd{\ps@pprintTitle}
  {Preprint to be submitted to}
  {Published in: }
  {}{}
\makeatother

\begin{document}
	
\begin{frontmatter}

\title{A Quantum Spectral Solver for Periodic Incompressible Stokes Flow}
\date{\today}

\author[inst1]{Juan M. Gimenez\corref{cor1}}
\ead{jmgimenez@cimne.upc.edu}
\cortext[cor1]{Corresponding author}

\author[inst2]{Fehmi Cirak}

\author[inst3,inst1]{Michael Ortiz}

\address[inst1]{Centre Internacional de Mètodes Numèrics a l'Enginyeria (CIMNE), Barcelona, 08014, Spain}
\address[inst2]{Department of Engineering, University of Cambridge, Cambridge, CB2 1PZ, UK}
\address[inst3]{Division of Engineering and Applied Science, California Institute of Technology, Pasadena, CA 91125, USA}

\journal{Comput. Methods Appl. Mech. Eng., 463:119453, 2027. DOI: 10.1016/j.cma.2026.119453}

\begin{abstract}
	We present a quantum spectral solver for the steady incompressible Stokes equations on a two-dimensional periodic domain. The method uses the Quantum Fourier Transform as a coherent change of basis and exploits the resulting spectral structure of the Stokes operator: the Laplacian becomes diagonal, while incompressibility is enforced mode by mode through a Helmholtz projection. In two dimensions, this projection is realized by a mode-dependent rotation from Cartesian velocity components to longitudinal--transverse coordinates, followed by component-conditioned inverse-Laplacian scaling. The velocity and pressure fields are encoded as quantum states over Fourier modes and physical components, and the corresponding spectral factors are implemented through polynomially encoded amplitude blocks.
	The construction extends recent quantum spectral methods in computational mechanics to an incompressible flow operator with explicit pressure--velocity splitting and divergence-free projection. The approach is also compatible with multiscale finite-element architectures in which quantum parallelism can simultaneously update all representative volume element (RVE) states. Numerical verification includes a steady vortex, a regularized periodic force-dipole benchmark, and an RVE-inspired Kolmogorov-like fluctuation benchmark. The latter illustrates how the circuit can recover a homogenized kinetic-energy observable without reconstructing the full velocity field, consistent with the role of averaged quantities in multiscale flow calculations. 
	Subject to application-dependent assumptions about efficient state preparation and observable estimation, 
	the circuit has polylogarithmic dependence on the grid resolution, with the polynomial degree and tile count appearing as explicit approximation and implementation parameters.
\end{abstract}

\begin{keyword}
	Quantum computing \sep Spectral methods \sep Stokes flow \sep Homogenization \sep Multiscale mechanics
\end{keyword}
\end{frontmatter}

\section{Introduction}

Periodic incompressible flow problems arise naturally in spectral discretizations, benchmark calculations, and representative-volume-element (RVE) models. In many such settings, the desired output is not only the resolved velocity and pressure field, but also averaged or homogenized quantities derived from that field, such as kinetic-energy measures, fluctuation stresses, fluxes, or coarse-scale response coefficients. This motivates formulations that exploit the structure of periodic incompressible operators while supporting both field reconstruction and observable estimation. More generally, periodic incompressible problems provide a useful setting in which operator structure, spectral representations, and homogenizable outputs can be studied in a controlled way. Related multiscale-flow methodologies likewise use representative incompressible-flow calculations to pass reduced or averaged information to a coarse-scale solver~\cite{Idelsohn2024PDNS,Gimenez2025PDNS}.

Quantum computing offers a route for rethinking such structured scientific-computing tasks. By representing an \(N\)-dimensional vector as amplitudes of a quantum state over \(n=\log_2 N\) qubits, the circuit model provides a compact state representation and a computational framework based on coherent unitary transformations and measurement~\cite{Nielsen2010}. For selected problem classes, quantum algorithms can manipulate high-dimensional states with resources that scale polynomially in \(n\), although any claim of advantage must be assessed for the complete workflow, including state preparation, oracle or block construction, circuit depth, success probability, observable-estimation cost, noise, and fault-tolerance overhead. Recent perspectives on scalable quantum computing identify block encodings and polynomial transformations as modular primitives that connect abstract quantum algorithms with analyzable resource estimates for computational-science applications~\cite{Joven2025}. \hl{The primary formulations underlying these primitives include the quantum linear-system algorithm, quantum singular value transformation, qubitization, and amplitude estimation and amplification~\cite{HHL2009,Gilyen2019,LowChuang2019,Brassard2002}.} This message is consistent with the emerging quantum computational-mechanics literature: useful quantum algorithms are unlikely to result from a literal gate-by-gate translation of classical solvers, and instead require operator factorizations that are naturally compatible with amplitude representations and unitary circuit constructions~\cite{Liu2024}. \hl{Recent studies also address distinct computational-mechanics bottlenecks through quantum-accelerated data-driven nearest-neighbour searches and decomposition-free variational linear solvers~\cite{Xu2024,Xu2025}.}

\hl{Spectral methods provide one such compatible structure in the established Fourier framework for incompressible flow~\cite{Canuto1988,Trefethen2000}.} On uniform grids, Fourier representations transform constant-coefficient differential operators into diagonal or block-diagonal operators in frequency space, while the change of basis between physical and spectral coefficients is implemented coherently by the Quantum Fourier Transform (QFT). Liu et al.~\cite{Liu2024} used this observation to construct a QFT-based quantum solver for periodic RVE problems in computational homogenization, combining spectral operator evaluation, piecewise polynomial approximation, and polynomial function encoding. Febrianto et al.~\cite{Febrianto2026} extended the same spectral philosophy to non-periodic boundary value problems by using antisymmetric extensions and quantum sine transforms to impose Dirichlet boundary conditions. Wang et al.~\cite{Wang2026} then introduced a quantum--classical multiscale finite-element framework in which the collection of microscopic RVE problems associated with macroscopic quadrature points is promoted to a quantum superposition. In that quantum accelerated multiscale finite element (QAFE$^2$) setting, quantum parallelism enables the simultaneous update of all quadrature points in the model, while entanglement preserves the association between each macroscopic quadrature state and its microscopic RVE response. These developments show that quantum spectral methods in computational mechanics can handle scalar spectral operators, non-periodic transforms, and concurrent RVE evaluations. They also clarify why efficient state preparation, block encoding, and polynomial amplitude loading are enabling components of quantum PDE solvers rather than secondary implementation details~\cite{GonzalezConde2025,Sanchez2023}. \hl{Their scope, however, does not include incompressible flow, for which enforcing a divergence-free velocity field is an additional operator-level requirement.}

Quantum algorithms for fluid mechanics have developed along several complementary lines. Succi et al.~\cite{Succi2023} frame this area as the use of quantum computers for classical fluid physics and review routes based on lattice-gas and lattice-Boltzmann formulations, nonlinear ODE solvers, variational methods, Carleman embeddings, and hydrodynamic reformulations inspired by Schr\"odinger dynamics. Tennie et al.~\cite{tennie2025quantum} provide a perspective focused on nonlinear differential equations and turbulence, discussing both algorithmic approaches to nonlinear dynamics and the hardware challenges that must be overcome before quantum computing can become relevant for large-scale turbulent-flow simulations. Raisuddin and De~\cite{Raisuddin2026} provide a complementary computational-mechanics review of quantum linear-system solvers, Hamiltonian simulation, ODE and PDE algorithms, block encoding, qubitization, quantum signal processing, amplitude amplification, and variational algorithms, emphasizing the gap between idealized algorithmic speedups and deployable engineering workflows. More specialized CFD contributions reflect this diversity: Carleman-linearized lattice-Boltzmann methods combined with HHL-type solvers have been explored for benchmark problems relevant to industrial CFD~\cite{Turro2025, luacuatucs2026surrogate}; hybrid quantum--classical schemes have used variational quantum routines for the pressure Poisson equation in incompressible Navier--Stokes projection methods~\cite{song25}; variational quantum algorithms have been proposed for potential and Stokes flows~\cite{LiuStokes2024}; and spectral quantum algorithms have been developed for passive-scalar advection--diffusion in prescribed shear flows with periodic, Neumann, and Dirichlet boundary conditions~\cite{Pfeffer2025}. 
More recent work has also focused on practical ingredients required for scalable quantum CFD implementations. Schalkers et al.~\cite{schalkers2024momentum} introduced a quantum momentum-exchange method for estimating forces on immersed objects without reconstructing the full flow field, while Georgescu et al.~\cite{georgescu2026efficient} developed a coherent strategy for imposing solid-boundary conditions on partitioned complex domains. \hl{Chen et al.~\cite{Chen2024} developed a hybrid Iterative QLS--subspace approach for near-term hardware, demonstrated for steady Poiseuille flow and unsteady acoustic-wave propagation. Other conceptually distinct approaches include tensor-programmable circuits for time-dependent differential equations, QSP/QSVT constructions for advection--diffusion and Koopman--von Neumann dynamics, and quantum-compatible homotopy linearization for nonlinear PDEs~\cite{Siegl2026,Novikau2025,Xue2025}.} Taken together, these efforts illustrate the different routes being developed for quantum algorithms for fluid mechanics.

A QFT-based spectral route for quantum CFD raises a specific operator-level question: how to represent incompressibility, pressure recovery, and non-unitary Stokes response within a structure-preserving quantum formulation. This route is attractive because, in the periodic setting, the Fourier transform diagonalizes the Laplacian and reduces incompressibility to a mode-wise Helmholtz projection. The steady periodic Stokes problem is therefore the minimal setting in which this structure appears in closed form, yielding a physically interpretable spectral operator that separates solenoidal velocity from pressure-generating modes and remains compatible with RVE-based multiscale architectures. \hl{Accordingly, the role of the present work is to establish and validate this incompressible spectral block as a prerequisite for, but not a full realization of, more general quantum CFD and multiscale solvers. The prospective multiscale challenge is the repeated evaluation of microscopic RVE responses associated with macroscopic quadrature points and the extraction of reduced quantities for the coarse-scale update.}

The objective of the present work is to develop, to the best of the authors' knowledge, the first QFT-based quantum spectral solver for the steady, periodic, two-dimensional incompressible Stokes equations. Starting from the classical Fourier reduction of the Stokes system, the algorithm encodes Fourier modes and velocity components in quantum registers, applies a shared multidimensional Fourier transform, realizes incompressibility through a mode-dependent longitudinal--transverse decomposition, and implements the velocity and pressure spectral operators through polynomially encoded amplitude blocks. 

The main contribution is threefold. First, we derive a quantum-compatible spectral factorization of the periodic incompressible Stokes equations that separates Helmholtz projection, inverse-Laplacian scaling, and pressure recovery. Second, we give an explicit circuit architecture with mode registers, a component register, a mode-dependent rotation, component-conditioned polynomial spectral blocks, and approximation parameters exposed through polynomial degree and tiling. Third, we validate the complete pipeline on a manufactured Taylor--Green check, a regularized periodic force-dipole benchmark, and an RVE-inspired Kolmogorov-like fluctuation benchmark, including recovery of a homogenized kinetic-energy observable without full-field reconstruction.

The remainder of the paper is organized as follows. Section~\ref{sec:classical-formulation} introduces the classical periodic Stokes problem and its spectral formulation. Section~\ref{sec:quantum-formulation} presents the quantum encoding, circuit construction, and implementation details. Section~\ref{sec:results} reports numerical results for the verification and RVE-inspired benchmarks. Finally, Section~\ref{sec:conclusions} summarizes the main findings and outlines directions for extending the proposed Stokes spectral block toward richer incompressible-flow and multiscale settings.

\section{Classical Spectral Formulation}
\label{sec:classical-formulation}

In this section, we recall the classical Fourier spectral formulation of the
two-dimensional steady Stokes equations on a periodic domain. The purpose is not
to introduce a new classical discretization, but to fix the notation and the
operator factorization used by the quantum construction in the following
sections. In Fourier space, the Stokes operator reduces mode by mode to a scalar
inverse-Laplacian operator composed with a mode-dependent incompressibility
projector.

\subsection{Stokes system and Fourier representation}

Let \(\Omega=[0,L]^2\) be a periodic square domain. The steady Stokes problem
for an incompressible Newtonian fluid reads
\begin{align}
	-\mu \Delta \bm u + \nabla p &= \bm f,
	\hspace{.5cm} \text{in } \Omega, \\
	\nabla \cdot \bm u &= 0,
	 \hspace{.5cm} \text{in } \Omega,
\end{align}
where \(\mu>0\) is the dynamic viscosity, \(\bm u=(u_0,u_1)^T\) is the velocity,
\(p\) is the pressure, and \(\bm f=(f_0,f_1)^T\) is a prescribed body force.
The pressure is defined up to an additive constant.

We expand all periodic fields in the same Fourier basis,
\begin{equation}
	g(\bm x)
	=
	\sum_{\bm k\in (2\pi/L)\mathbb Z^2}
	\hat g(\bm k) e^{i\bm k\cdot \bm x},
	\qquad
	g\in\{\bm u,p,\bm f\}.
\end{equation}
Here \((2\pi/L)\mathbb Z^2\) denotes the physical Fourier lattice of the periodic domain.
Using
\begin{equation}
	\nabla \mapsto i\bm k,
	\qquad
	\Delta \mapsto -|\bm k|^2,
	\qquad
	\nabla\cdot \mapsto i\bm k\cdot,
\end{equation}
the Stokes equations decouple mode by mode into
\begin{align}
	\mu |\bm k|^2 \hat{\bm u}(\bm k)
	+ i\bm k \hat p(\bm k)
	&=
	\hat{\bm f}(\bm k),
	\label{eq:fourier-stokes-momentum}
	\\
	\bm k \cdot \hat{\bm u}(\bm k)
	&=
	0.
	\label{eq:fourier-stokes-incompressibility}
\end{align}

\subsection{Pressure elimination and Helmholtz projection}

For every nonzero mode \(\bm k\neq \bm 0\), taking the scalar product of
\eqref{eq:fourier-stokes-momentum} with \(\bm k\) and using
\eqref{eq:fourier-stokes-incompressibility} gives
\begin{equation}
	i|\bm k|^2 \hat p(\bm k)
	=
	\bm k\cdot \hat{\bm f}(\bm k).
\end{equation}
Therefore,
\begin{equation}
	\hat p(\bm k)
	=
	-i
	\frac{\bm k\cdot \hat{\bm f}(\bm k)}
	{|\bm k|^2}.
	\label{eq:pressure-fourier-solution}
\end{equation}
Substitution into the momentum equation~\eqref{eq:fourier-stokes-momentum} yields
\begin{equation}
	\hat{\bm u}(\bm k)
	=
	\frac{1}{\mu |\bm k|^2}
	\left(
	I
	-
	\frac{\bm k\bm k^T}{|\bm k|^2}
	\right)
	\hat{\bm f}(\bm k),
	\qquad
	\bm k\neq \bm 0.
	\label{eq:stokes-fourier-projector}
\end{equation}
The matrix
\begin{equation}
	P_\perp(\bm k)
	=
	I
	-
	\frac{\bm k\bm k^T}{|\bm k|^2}
	\label{eq:helmholtz-projector}
\end{equation}
is the Fourier representation of the Helmholtz projector onto the subspace
orthogonal to \(\bm k\). Thus, the velocity is obtained by first removing the longitudinal component of
the force and then applying the scalar inverse Laplacian,
\begin{equation}
	\hat{\bm u}(\bm k)
	=
	\frac{1}{\mu |\bm k|^2}
	P_\perp(\bm k)
	\hat{\bm f}(\bm k).
\end{equation}
For the zero mode, the inverse factor \(1/|\bm k|^2\) is singular. The periodic
Stokes problem is solvable, modulo the pressure gauge, only when the forcing
satisfies the usual zero-mean compatibility condition. We therefore impose the
zero-mode convention
\begin{equation}
	\hat{\bm u}(\bm 0)=\bm 0,
	\qquad
	\hat p(\bm 0)=0,
\end{equation}
with the latter fixing the arbitrary additive constant in the pressure.

\subsection{Discrete spectral coordinates}\label{sec:discrete-spectral}

For an \(N\times N\) periodic grid with \(N=2^n\), we label the discrete Fourier
frequencies using the unsigned indices \(k^0,k^1\in\{0,\ldots,N-1\}\) and the
signed integer map
\begin{equation}
	r(k)
	=
	\begin{cases}
		k, & 0\le k < N/2, \\
		k-N, & N/2 \le k < N.
	\end{cases}
	\label{eq:r-map}
\end{equation}
The wavevector associated with the grid index pair \((k^0,k^1)\) is
\begin{equation}
	\bm k(k^0,k^1)
	=
	\frac{2\pi}{L}
	\begin{pmatrix}
		r(k^0) \\
		r(k^1)
	\end{pmatrix}.
\end{equation}
In the following, when no ambiguity arises, we write \(\bm k\equiv \bm k(k^0,k^1)\), therefore
\begin{equation}
	|\bm k|^2
	=
	\frac{4\pi^2}{L^2}
	\left[
	r(k^0)^2+r(k^1)^2
	\right].
\end{equation}
The numerical method uses the corresponding band-limited approximation
\begin{equation}
	\bm u_N(\bm x)
	=
	\sum_{k^0=0}^{N-1}\sum_{k^1=0}^{N-1}
	\hat{\bm u}\!\left(\bm k \right)
	\exp\!\left(i\bm k \cdot \bm x\right),
	\label{eq:band-limited-velocity}
\end{equation}
with analogous expansions for \(p_N\) and \(\bm f_N\) and $\bm x = (x_0, x_1)^T$. The subscript \(N\) denotes the discrete spectral truncation only; it does not
represent a different continuous Stokes problem. \hl{The mode-wise Stokes relations are independent of the global normalization used for the discrete Fourier coefficients.}

The scalar velocity Green operator used below is then
\begin{equation}
	\Gamma(k^0,k^1)
	=
	\frac{1}{\mu |\bm k|^2}
	=
	\frac{L^2}
	{4\pi^2 \mu \left[r(k^0)^2+r(k^1)^2\right]},
	\qquad
	(k^0,k^1)\neq (0,0),
\end{equation}
with
\begin{equation}
	\Gamma(0,0)=0.
\end{equation}
Thus, the discrete spectral Stokes solver applies, for each nonzero Fourier
mode, the matrix
\begin{equation}
	\Gamma(k^0,k^1)\, P_\perp(\bm k),
\end{equation}
or equivalently its rotation-based factorization in
\eqref{eq:stokes-rotation-form}.

\subsection{Rotation-based representation of the projector}

In two dimensions, the Helmholtz projector can also be expressed through a
mode-dependent rotation. For each nonzero mode, define the orthonormal directions
\begin{equation}
	\bm e_\parallel(\bm k)
	=
	\frac{1}{\rho(k^0,k^1)}
	\begin{pmatrix}
		r(k^0) \\ r(k^1)
	\end{pmatrix},
	\qquad
	\bm e_\perp(\bm k)
	=
	\frac{1}{\rho(k^0,k^1)}
	\begin{pmatrix}
		-r(k^1) \\ r(k^0)
	\end{pmatrix},
\end{equation}
where $\rho(k^0,k^1) = \sqrt{r(k^0)^2+r(k^1)^2}$.
The matrix
\begin{equation}
	R(\bm k)
	=
	\begin{pmatrix}
		\bm e_\parallel(\bm k)^T \\
		\bm e_\perp(\bm k)^T
	\end{pmatrix}
	=
	\frac{1}{\rho(k^0,k^1)}
	\begin{pmatrix}
		r(k^0) & r(k^1) \\
		-r(k^1) & r(k^0)
	\end{pmatrix}
\end{equation}
rotates vector components from the Cartesian basis to the
\((\parallel,\perp)\) basis associated with the Fourier mode \(\bm k\). In this
basis, the incompressibility constraint removes the longitudinal component and
keeps only the transverse one. Hence,
\begin{equation}
	P_\perp(\bm k)
	=
	R(\bm k)^T
	\begin{pmatrix}
		0 & 0 \\
		0 & 1
	\end{pmatrix}
	R(\bm k).
	\label{eq:projector-rotation-form}
\end{equation}
Consequently, the Fourier-space Stokes solution may be written as
\begin{equation}
	\hat{\bm u}(\bm k)
	=
	R(\bm k)^T
	\begin{pmatrix}
		0 & 0 \\
		0 & \dfrac{1}{\mu |\bm k|^2}
	\end{pmatrix}
	R(\bm k)
	\hat{\bm f}(\bm k),
	\qquad
	\bm k\neq \bm 0.
	\label{eq:stokes-rotation-form}
\end{equation}
This factorization then consists of a
rotation into the longitudinal--transverse basis, a transverse component-conditioned inverse-Laplacian scaling, and a rotation back to Cartesian
coordinates.

Equivalently, introducing
\begin{equation}
	e_0(\bm k)=\frac{r(k^0)}{\rho(k^0,k^1)},
	\qquad
	e_1(\bm k)=\frac{r(k^1)}{\rho(k^0,k^1)},
\end{equation}
or the polar angle \(\theta(\bm k)\) such that
\begin{equation}
	\cos \theta(\bm k)=e_0(\bm k),
	\qquad
	\sin \theta(\bm k)=e_1(\bm k),
\end{equation}
one obtains
\begin{equation}
	R(\bm k)
	=
	\begin{pmatrix}
		\cos\theta(\bm k) & \sin\theta(\bm k) \\
		-\sin\theta(\bm k) & \cos\theta(\bm k)
	\end{pmatrix}.
\end{equation}
This angular representation is particularly useful for the quantum circuit
implementation, where the projector can be realized through mode-dependent
rotations and a diagonal transverse-mode scaling.

\section{Quantum Spectral Formulation and Circuit Realization}
\label{sec:quantum-formulation}

We now describe the quantum implementation of the Fourier spectral Stokes
solver introduced in Section~\ref{sec:classical-formulation}. The classical
algorithm starts from a periodic forcing field, transforms it to Fourier space,
applies mode-wise Stokes operators, and transforms the selected result back to
physical space. The quantum construction follows the same sequence, but fields
are represented by amplitudes of a normalized state and every intermediate map
must be embedded in a unitary operation. Consequently, non-unitary operations
such as division by \(|\bm k|^2\), post-selection, and observable extraction
are implemented by using ancillary registers, postselecting the corresponding
success subspaces, and applying explicit rescaling factors.

Throughout this section, a non-unitary spectral map should be read as
being embedded in a larger unitary acting on the physical registers and on
ancillas. The desired Stokes response is the component of the final state in a
postselected subspace of the ancillary registers; the physical normalization is then restored using the norm of the input forcing and the scale factors used in the amplitude encodings. This convention is standard in block-encoding and
polynomial-function-loading constructions, and is used here only to make the
Stokes spectral factors compatible with quantum computing~\cite{Joven2025,GonzalezConde2025}.

The implementation builds on the pattern used in the quantum spectral solvers in computational mechanics
of Liu et al.~\cite{Liu2024} and Febrianto et al.~\cite{Febrianto2026}:
state preparation encodes the input data, the quantum Fourier transform replaces
the FFT as a coherent change of basis, and mode-dependent spectral multipliers
are approximated by polynomial encodings. The Stokes-specific addition is the
rotation-based representation of the incompressibility projector. After the
Fourier transform, a mode-dependent component rotation separates longitudinal
and transverse amplitudes; the velocity, pressure, and homogenized observables
are then obtained from different postselection of the same prepared
spectral state. 

\subsection{Quantum registers, input state, and normalization}

We use Dirac notation to identify the basis states and amplitudes carried
by the quantum register. A basis state is an element of the computational basis;
the coefficient multiplying that basis state is the amplitude that stores a
field value after normalization. For the two-dimensional grid, the circuit uses two \(n\)-qubit grid/mode
registers, denoted by \(k^0\) and \(k^1\). 
The basis states are denoted as
\begin{equation}
|k^0,k^1\rangle := |k^0\rangle |k^1\rangle,
\qquad
0\le k^0,k^1<N .
\end{equation}
Each \hl{grid/mode register} is represented by \(n\) qubits.
For \(N=2^n\),
\begin{equation}
k^a=\sum_{\ell=0}^{n-1}2^\ell k^a_\ell,
\qquad a=0,1,
\end{equation}
so that
\begin{equation}
|k^0,k^1\rangle
=
|k^0\rangle|k^1\rangle
=
|k^0_0k^0_1\cdots k^0_{n-1}\rangle
|k^1_0k^1_1\cdots k^1_{n-1}\rangle .
\end{equation}
Here \(k^0\) and \(k^1\) denote the two registers, while \(k^a_\ell\) denotes
the state of the \(\ell\)th qubit of register \(k^a\). The same labels \(k^0,k^1\) are used
before and after the QFT. Before applying \(\mathrm{QFT}_{2D}\), they label physical grid
points at which the forcing is sampled. After applying \(\mathrm{QFT}_{2D}\), they label
discrete Fourier modes, with signed coordinates obtained through the
relabelling map \(r(k^a)\) introduced in Section~\ref{sec:classical-formulation}.
Thus the register names and basis labels are fixed throughout the circuit; only
their interpretation changes across the QFT. No new wavevector notation is
introduced in this section: \(\bm k\) denotes the Fourier mode associated with
the current labels when the state is in spectral form.

The complete register is written as
\begin{equation}
|k^0\rangle|k^1\rangle|c\rangle|t\rangle|w\rangle,
\end{equation}
where the component qubit stores the two Cartesian components as
\begin{equation}
	|0\rangle_c \equiv 0\text{-component},
	\qquad
	|1\rangle_c \equiv 1\text{-component}.
\end{equation}
	
The target qubit is used for amplitude encoding of scalar spectral factors,
while the work register stores ancillary tile labels, polynomial controls, and temporary arithmetic. 
\hl{It contains $n+1$ qubits: two reusable tile-membership flags and $n-1$ ancillary qubits for reversible comparisons. Thus, for $N=2^n$, the circuit uses $2n$ grid/mode qubits, one component qubit, one target qubit, and $n+1$ work qubits, for a total of $3n+3$ qubits.}

The initialized register reads
\begin{equation}
|0\rangle_{k^0}|0\rangle_{k^1}|0\rangle_c|0\rangle_t|0\rangle_w .
\end{equation}
All operators below are unitary on the entire register. The physical Stokes
operator is not unitary by itself, so the desired field is recovered from a
specified branch of the final state and then rescaled. This is the same
principle used for polynomial function encodings in previous quantum spectral methods in computational mechanics
constructions, but here the branch also identifies the vector
component selected by the incompressibility projection.

The state-preparation unitary \(U_I\) loads the normalized physical-space
forcing field,
\begin{equation}
	U_I |0\rangle_{k^0}|0\rangle_{k^1}|0\rangle_c|0\rangle_t|0\rangle_w
	=
	|\psi_f\rangle |0\rangle_t|0\rangle_w,
\end{equation}
with
\begin{equation}
	|\psi_f\rangle
	=
	\frac{1}{\|\bm f\|}
	\sum_{k^0,k^1}
	\left[
	f_0(k^0,k^1)|k^0\rangle|k^1\rangle|0\rangle_c
	+
	f_1(k^0,k^1)|k^0\rangle|k^1\rangle|1\rangle_c
	\right].
	\label{eq:forcing-state}
\end{equation}
In this input state, \(f_i(k^0,k^1)\) denotes the value of the \(i\)th Cartesian
component of the forcing sampled at the physical grid point labelled by
\((k^0,k^1)\); it is not a Fourier coefficient. Here
\begin{equation}
	\|\bm f\|^2
	=
	\sum_{k^0,k^1}
	\left(
	|f_0(k^0,k^1)|^2
	+
	|f_1(k^0,k^1)|^2
	\right)
\end{equation}
is the Euclidean norm of the discrete forcing vector. 
\hl{This norm depends on the grid resolution. On the unit square, $F^2=\|\bm f\|^2/N^2$ is the physical-grid quadrature approximation of the continuum squared forcing norm.}
Thus $U_I$ prepares a
unit-norm quantum state proportional to the physical forcing. If the physical
forcing is not already unit-normalized, recovered fields and observables must
be rescaled by the appropriate powers of \(\|\bm f\|\) and by the
amplitude-encoding scale factors introduced below. 
The component qubit is a structural part of the discretized Stokes operator: it represents the local vector degree of freedom and allows both Cartesian
components to share the same spatial and spectral registers.
\hl{In practice, $U_I$ can be realized by controlled state-preparation unitaries~\cite{PleschBrukner2011,Zhang2022,Sanchez2023}. For forcing components whose grid samples admit a compact piecewise-polynomial approximation, Liu et al.~\cite{Liu2024} provide a concrete function-encoding route in which controlled rotations load the approximate input amplitudes from a uniform superposition over grid labels, thereby avoiding the explicit loading of $O(2^{2n})$ unrelated amplitudes.}

\subsection{Polynomial encoding, unitary embedding, and tiled implementation}
\label{subsec:poly}
Several unitary blocks that will be introduced below depend on the Fourier label stored in the spectral
registers. In a classical spectral code, functions such as \(1/|\bm k|^2\) or
angular projection coefficients are evaluated directly on the FFT grid. In a
quantum implementation, the same functions must be evaluated coherently on a
superposition of mode labels. Following the polynomial function-encoding
viewpoint used by Liu et al.~\cite{Liu2024} and the piecewise spectral
construction of Febrianto et al.~\cite{Febrianto2026}, we approximate these
mode-dependent maps by low-degree polynomials on a tiling of the spectral plane.

A tile \(\tau\) is a rectangular subset of the discrete index plane
\((k^0,k^1)\). On each tile, the indices are mapped affinely to local
coordinates \(\xi_0,\xi_1\in[-1,1]\). Let \(a\) denote a mode-dependent function
to be encoded. On \(\tau\), we use the tensor-product polynomial convention
\begin{equation}
	a(\bm k)
	\approx
	a_\tau(\bm k)
	=
	\sum_{\alpha=0}^{p}
	\sum_{\beta=0}^{p}
	c_{\alpha\beta}^{(\tau)}
	\xi_0(\bm k)^\alpha
	\xi_1(\bm k)^\beta,
	\label{eq:tile-polynomial}
\end{equation}
where the same degree \(p\) is used in both directions. The coefficients may be
obtained in a tensor-product Chebyshev basis and then converted to monomials for
circuit synthesis. This is the convention used in the numerical examples.

The tiled implementation follows the same logic as the piecewise polynomial
encoding in Liu et al.~\cite{Liu2024}. Comparator subcircuits tag whether the
current spectral label belongs to a tile, the corresponding local polynomial is applied conditionally on that tag, and the tile flags and temporary arithmetic are then uncomputed so that the ancillary work register returns to \(|0\rangle_w\) on the selected branch. Tiling is therefore both an approximation choice and a circuit-design choice. More tiles improve local approximation of nonsmooth or rapidly varying spectral operators, but each tile introduces comparator, control, and uncomputation overhead. Likewise, increasing \(p\) adds controlled rotations associated with additional monomial terms; after transpilation this appears as growth in both CNOT and \(U_3\) counts, consistent with the gate-count scaling discussed for polynomial encodings in Liu et al.~\cite{Liu2024}.

For scalar spectral operators such as \(\Gamma\) and \(\Lambda\) \hl{further introduced in Eqs.~\eqref{eq:gamma} and~\eqref{eq:lambda}}, the
polynomial is encoded as the amplitude of a target qubit.
We use the standard single-qubit Pauli-$Y$ rotation
\begin{equation}
R_Y(\varphi)
=
\exp\left(-i\frac{\varphi}{2}Y\right)
=
\begin{pmatrix}
	\cos(\varphi/2) & -\sin(\varphi/2) \\
	\sin(\varphi/2) & \cos(\varphi/2)
\end{pmatrix}.
\end{equation}
In the scalar blocks below, this rotation acts on the target qubit $t$.
Later, the same gate family is used on the component qubit $c$ to implement
the longitudinal--transverse change of basis. Therefore, the elementary
rotation is
\begin{equation}
	R_Y(2\alpha_a(\bm k))|0\rangle_t
	=
	\cos\alpha_a(\bm k)|0\rangle_t
	+
	\sin\alpha_a(\bm k)|1\rangle_t .
	\label{eq:elementary-amplitude-rotation}
\end{equation}

For amplitude encoding, the \(|1\rangle\) branch is linearized in the target function \hl{
by assuming a small-angle regime achieved via a small scale parameter $\varepsilon_a$. Under this assumption, the sine function can be expressed as}
\begin{equation}
	\sin\alpha_a(\bm k)
	\simeq
	\alpha_a(\bm k)
	=
	\varepsilon_a a(\bm k),
	\qquad
	|\varepsilon_a a(\bm k)|\ll 1,
\label{eq:amplitude-linearization}
\end{equation}
\hl{where decreasing $\varepsilon_a$ reduces the
linearization error but also decreases the amplitude of the selected target
branch.}

With this linearization, the generic scalar unitary has the following
postselected success-branch action
\begin{equation}
	U_a:
	|k^0\rangle|k^1\rangle |c\rangle |0\rangle_t |0\rangle_w
	\mapsto
	|k^0\rangle|k^1\rangle |c\rangle
	\left[
	\sqrt{1-\varepsilon_a^2 a(\bm k)^2}|0\rangle_t
	+
	\varepsilon_a a(\bm k)|1\rangle_t
	\right]
	|0\rangle_w
	+
	Junk,
	\label{eq:generic-polynomial-encoding}
\end{equation}
where \(\varepsilon_a\) is chosen such that
\(|\varepsilon_a a(\bm k)|\le 1\) over the represented spectral domain.
In the actual synthesis the polynomial approximant is used for \(a\), and the
linearity of the sine is an approximation controlled by the selected scale and
polynomial accuracy. The work register returns to \(|0\rangle_w\) only on the
postselected success branch used for solution extraction. Here and below, $Junk$ denotes components in branches that are orthogonal to the selected solution branch, including branches with nonzero work registers or unsuccessful target-register outcomes.

\subsection{Shared spectral and component-space operators}

\subsubsection{Quantum Fourier transform}

The Fourier series in Section~\ref{sec:classical-formulation} uses the basis
\(e^{i\bm k\cdot\bm x}\). We therefore use the forward spectral transform with
negative phase,
\begin{equation}
	\mathrm{QFT}_{2D}
	=
	\mathrm{QFT}_{0,N}\otimes \mathrm{QFT}_{1,N},
\end{equation}
where
\begin{equation}
	\mathrm{QFT}_N |\ell\rangle
	=
	\frac{1}{\sqrt{N}}
	\sum_{m=0}^{N-1}
	e^{-2\pi i \ell m/N}|m\rangle ,
	\qquad \ell = 0,\ldots,N-1 .
	\label{eq:qft-sign-convention}
\end{equation}
This unitary is the adjoint of the positive-phase convention commonly used for
the textbook QFT. The choice in \eqref{eq:qft-sign-convention} is purely
for convenience; it makes the amplitudes after \(\mathrm{QFT}_{2D}\) coincide with the unitary discrete Fourier coefficients of the expansion with \(e^{+i\bm k\cdot\bm x}\), and \(\mathrm{QFT}_{2D}^\dagger\) is the corresponding inverse transform.

\hl{For an $N\times N$ grid, the amplitudes produced by $\mathrm{QFT}_{2D}$ use the unitary normalization $N^{-1}$, whereas the Fourier-series coefficients in Section~\ref{sec:discrete-spectral} use the conventional $N^{-2}$ normalization. Thus, the unitary coefficients are larger by a factor $N$; the inverse QFT supplies the reciprocal factor when reconstructing physical-grid values.}

Here \(\ell\) and \(m\) are generic input and output labels of the one-dimensional transform. In the Stokes circuit, these labels correspond to the same register labels \(k^a\) before and after the transform, with their interpretation changing from grid labels to Fourier-mode labels. Applied to the two grid/mode registers, the transform maps the physical-space forcing state to 
\begin{equation}
	(\mathrm{QFT}_{2D}\otimes I_c)|\psi_f\rangle
	=
	\frac{1}{\|\bm f\|}
	\sum_{k^0,k^1}
	\left[
	\hat f_0(k^0,k^1)|k^0\rangle|k^1\rangle|0\rangle_c
	+
	\hat f_1(k^0,k^1)|k^0\rangle|k^1\rangle|1\rangle_c
	\right].
	\label{eq:qft-action}
\end{equation}
where the hats denote unitary-DFT coefficients consistent with \eqref{eq:qft-sign-convention}. After \(\mathrm{QFT}_{2D}\), the same labels \((k^0,k^1)\) are interpreted as Fourier-mode labels and mapped to signed wavevectors by the
relabelling function \(r\) defined in~\eqref{eq:r-map}.

The QFT is the point at which the quantum algorithm inherits the structure of a spectral method. It replaces the classical FFT stage by a coherent basis change on the \(2n\) qubits that label the \(N^2\) grid points. The transform acts on all amplitudes of the encoded forcing state at once, and the cost of the standard circuit grows polynomially in \(n=\log_2 N\), rather than with the number of grid points. After the transform, differential operators become mode-dependent algebraic multipliers.
Because the two Cartesian components are stored as amplitudes of the same
quantum state, with the component qubit distinguishing them, the single
operation \(\mathrm{QFT}_{2D}\otimes I_c\) transforms both components simultaneously.
The QFT itself is a standard and extensively studied quantum algorithm
\cite{Nielsen2010}. In this work it is therefore treated as an established
basis-change primitive rather than as a new circuit contribution; practical
software stacks such as Qiskit provide QFT implementations that can be used or
specialized inside larger algorithms~\cite{Qiskit2024}.

\subsubsection{Scalar spectral blocks}

The velocity computation uses the Stokes Green operator
\begin{equation}
	\Gamma(\bm k)
	=
	\frac{1}{\mu |\bm k|^2},
	\qquad
	\Gamma(\bm 0)=0,
    \label{eq:gamma}
\end{equation}
while the pressure computation uses
\begin{equation}
	\Lambda(\bm k)
	=
	\frac{1}{|\bm k|},
	\qquad
	\Lambda(\bm 0)=0.
    \label{eq:lambda}
\end{equation}
Both functions are diagonal in the spectral registers and are implemented with
the polynomial amplitude-encoding convention in
\eqref{eq:generic-polynomial-encoding}, analogously to the scalar spectral
multipliers used for the Poisson operator in Febrianto et al.~\cite{Febrianto2026}.
We denote the corresponding unitaries by \(U_\Gamma\) and \(U_\Lambda\). On the
selected branch,
\begin{align}
	U_\Gamma:\;
	|k^0\rangle|k^1\rangle|1\rangle_c|0\rangle_t|0\rangle_w
	&\mapsto
	\varepsilon_\Gamma\Gamma(\bm k)
	|k^0\rangle|k^1\rangle|1\rangle_c|1\rangle_t|0\rangle_w + Junk, \\
	U_\Lambda:\;
	|k^0\rangle|k^1\rangle|0\rangle_c|0\rangle_t|0\rangle_w
	&\mapsto
	\varepsilon_\Lambda\Lambda(\bm k)
	|k^0\rangle|k^1\rangle|0\rangle_c|1\rangle_t|0\rangle_w + Junk.
\end{align}
The first action is used on the transverse branch for velocity, and the second
on the longitudinal branch for pressure.

\subsubsection{\texorpdfstring{Angular component rotation \(U_R\)}{Angular component rotation UR}}

For each nonzero register label \((k^0,k^1)\), let
\(\bm k=\bm k(k^0,k^1)\) denote the associated Fourier wavevector, and define
\begin{equation}
	e_0(\bm k)
	=
	\frac{r(k^0)}{\sqrt{r(k^0)^2+r(k^1)^2}},
	\qquad
	e_1(\bm k)
	=
	\frac{r(k^1)}{\sqrt{r(k^0)^2+r(k^1)^2}}.
\end{equation}
The classical rotation from Cartesian components to
longitudinal--transverse components is given by
\begin{equation}
	R(\bm k)
	=
	\begin{pmatrix}
		e_0(\bm k) & e_1(\bm k)\\
		-e_1(\bm k) & e_0(\bm k)
	\end{pmatrix}.
	\label{eq:Rk-definition}
\end{equation}
Introducing the polar angle
\begin{equation}
	\theta(\bm k)
	=
	\operatorname{atan2}\!\left(r(k^1),r(k^0)\right),
\end{equation}
one has
\[
e_0(\bm k)=\cos\theta(\bm k),
\qquad
e_1(\bm k)=\sin\theta(\bm k),
\]
and therefore
\begin{equation}
	R(\bm k)
	=
	\begin{pmatrix}
		\cos\theta(\bm k) & \sin\theta(\bm k)\\
		-\sin\theta(\bm k) & \cos\theta(\bm k)
	\end{pmatrix}.
\end{equation}

Using the Pauli-\(Y\) rotation convention introduced in
Section~\ref{subsec:poly}, this component rotation is obtained by choosing
\begin{equation}
	\alpha_R(\bm k)
	=
	-2\theta(\bm k)
	=
	-2\operatorname{atan2}\!\left(r(k^1),r(k^0)\right).
	\label{eq:rotation-angle-alpha}
\end{equation}
Indeed,
\begin{equation}
	R_Y(\alpha_R(\bm k))
	=
	\begin{pmatrix}
		\cos\theta(\bm k) & \sin\theta(\bm k)\\
		-\sin\theta(\bm k) & \cos\theta(\bm k)
	\end{pmatrix}
	=
	R(\bm k).
	\label{eq:Ry-equals-Rk}
\end{equation}
For the zero mode, we set \(U_R=I_c\). This choice is immaterial for the
Stokes response, since the zero mode is assigned zero response by the scalar
spectral blocks.
Thus, applying \(U_R\) to the component register converts the Cartesian spectral
amplitudes into longitudinal and transverse amplitudes:
\begin{equation}
	U_R:
	|k^0,k^1\rangle
	\left(
	\hat f_0(\bm k)|0\rangle_c
	+
	\hat f_1(\bm k)|1\rangle_c
	\right) |0\rangle_t |0\rangle_w
	\mapsto
	|k^0,k^1\rangle
	\left(
	\hat f_\parallel(\bm k)|0\rangle_c
	+
	\hat f_\perp(\bm k)|1\rangle_c
	\right)  |0\rangle_t |0\rangle_w,
	\label{eq:UR-action}
\end{equation}
where
\begin{equation}
	\begin{pmatrix}
		\hat f_\parallel(\bm k)\\
		\hat f_\perp(\bm k)
	\end{pmatrix}
	=
	R(\bm k)
	\begin{pmatrix}
		\hat f_0(\bm k)\\
		\hat f_1(\bm k)
	\end{pmatrix}
	=
	R_Y(\alpha_R(\bm k))
	\begin{pmatrix}
		\hat f_0(\bm k)\\
		\hat f_1(\bm k)
	\end{pmatrix}.
	\label{eq:UR-component-action}
\end{equation}
Equivalently, after input preparation and the two-dimensional QFT,
\begin{equation}
	U_R(\mathrm{QFT}_{2D}\otimes I_c)U_I
	|0\rangle_{k^0}|0\rangle_{k^1}|0\rangle_c|0\rangle_t|0\rangle_w
	=
	\frac{1}{\|\bm f\|}
	\sum_{k^0,k^1}
	|k^0,k^1\rangle
	\left[
	\hat f_\parallel(\bm k)|0\rangle_c
	+
	\hat f_\perp(\bm k)|1\rangle_c
	\right]
	|0\rangle_t|0\rangle_w .
	\label{eq:shared-rotated-forcing}
\end{equation}
The amplitudes \(\hat f_\parallel\) and \(\hat f_\perp\) are the longitudinal branch \(c=0\) and the
transverse branch \(c=1\) generated by \(U_R\) from the Cartesian spectral
forcing. The zero mode is assigned zero response and is excluded from the
angular encoding.

For implementation, the target function for this block is not \(e_0\), \(e_1\),
\(\sin\theta\), or \(\cos\theta\), but the angle
\(\alpha_R(\bm k)\) defined in~\eqref{eq:rotation-angle-alpha}. On each tile
\(\tau\), this angle is approximated using the polynomial encoding framework of
Section~\ref{subsec:poly},
\begin{equation}
	\alpha_R(\bm k)
	\approx
	\alpha_{R,\tau}(\bm k).
	\label{eq:rotation-angle-polynomial}
\end{equation}
The implemented gate on tile \(\tau\) is therefore
\(R_Y(\alpha_{R,\tau}(\bm k))\), which approximates the ideal
mode-dependent rotation \(\hll{R_Y}(\alpha_R(\bm k))=R(\bm k)\).

The angular map has quadrant discontinuities and rapid low-mode variation, so
its tile layout is more sensitive than that of a smooth scalar Green operator.
\hl{Accordingly, the angular map is represented by quadrant-aware piecewise polynomials, with finer tiles allocated to the low-wavenumber region; the resulting layout is examined in Section~\ref{sec:results}.}
Increasing \(N\) mainly extends the high-frequency tail of this map. This is an
important distinction from the scalar amplitude blocks: no small-angle
approximation such as \(\sin\alpha\simeq \alpha\) is used for \(U_R\). The
polynomial approximates the angle, and the quantum gate applies the
corresponding trigonometric rotation exactly at the gate level.

\subsection{Velocity circuit}
\label{sec:velocity-circuit}

\begin{figure}[h]
	\centering
	\[
	\Qcircuit @C=1.0em @R=1.1em {
		\lstick{\lvert k^0 \rangle} & \multigate{2}{U_I} & \gate{\mathrm{QFT}_{0,N}} & \multigate{4}{U_R}  & \multigate{4}{U_\Gamma} & \multigate{4}{U_R^\dagger} & \gate{\mathrm{QFT}_{0,N}^\dagger} & \qw \\
		\lstick{\lvert k^1 \rangle} & \ghost{U_I}        & \gate{\mathrm{QFT}_{1,N}} & \ghost{U_R}                & \ghost{U_\Gamma}        & \ghost{U_R^\dagger}        & \gate{\mathrm{QFT}_{1,N}^\dagger} & \qw \\
		\lstick{\lvert c \rangle}   & \ghost{U_I}        & \qw            & \ghost{U_R}                & \ghost{U_\Gamma}        & \ghost{U_R^\dagger}        & \qw                    & \qw \\
		\lstick{\lvert t \rangle}   & \qw                & \qw            & \ghost{U_R}               & \ghost{U_\Gamma}        & \ghost{U_R^\dagger}        & \qw                    & \qw \\
		\lstick{\lvert w \rangle}   & \qw                & \qw            & \ghost{U_R}                & \ghost{U_\Gamma}        & \ghost{U_R^\dagger}        & \qw                    & \qw
	}
	\]
		\caption{Quantum circuit for the velocity field on registers
			\(k^0,k^1,c,t,w\). The shared preprocessing \(U_I\), \(\mathrm{QFT}_{2D}\), and \(U_R\)
			is followed by the component-conditioned scalar scaling
			\(U_\Gamma\) on the transverse branch \(c=1\). Then \(U_R^\dagger\)
			returns the component register to the Cartesian basis. The selected output
			branch is \(t=1\) and \(w=0\).}
	\label{fig:velocity-circuit}
\end{figure}
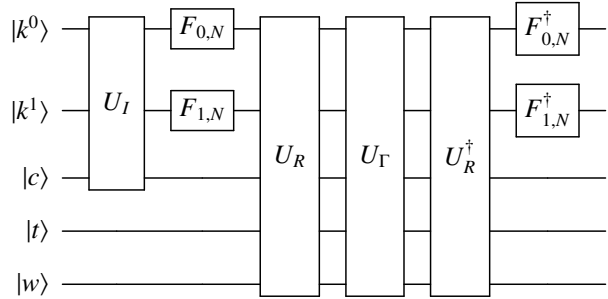

From the state presented in Eq.~\eqref{eq:shared-rotated-forcing}, the velocity circuit follows the transverse branch \(c=1\):
\(U_\Gamma\) acts component-conditionally on that branch, \(U_R^\dagger\)
returns the component register to Cartesian coordinates, and \(\mathrm{QFT}_{2D}^{\dagger}\)
reconstructs the physical field, as shown in
Figure~\ref{fig:velocity-circuit}. Its postselected operator is
\begin{equation}
	\mathcal U_{\bm u}
	=
	(\mathrm{QFT}_{2D}^{\dagger}\otimes I_c)
	U_R^\dagger
	U_\Gamma
	U_R
	(\mathrm{QFT}_{2D}\otimes I_c)
	U_I .
	\label{eq:velocity-operator}
\end{equation}

\hl{Before $U_R^\dagger$, $U_\Gamma$ acts on the transverse branch $c=1$; after $U_R^\dagger$, $c$ again labels Cartesian components. On the postselected ancillary-success branch $t=1$ and $w=0$, the final state after $\mathrm{QFT}_{2D}^{\dagger}$ is proportional to the physical-space velocity field,}
\begin{equation}
	|\psi_{\bm u}\rangle
	\propto
	\sum_{k^0,k^1}
	\left[
	u_0(k^0,k^1)
	|k^0\rangle|k^1\rangle|0\rangle_c
	+
	u_1(k^0,k^1)
	|k^0\rangle|k^1\rangle|1\rangle_c
	\right].
	\label{eq:velocity-final-state}
\end{equation}
The proportionality \hl{means that} the selected amplitudes include the known factor
\(\varepsilon_\Gamma/\|\bm f\|\). Then, these amplitudes should be rescaled by the inverse of this factor to obtain the physical values for velocity.

When \(U_R\) and \(U_\Gamma\) use the same tile partition, the tile membership
flags can be computed once and reused across \(U_R\), \(U_\Gamma\), and
\(U_R^\dagger\). This leaves the selected-branch operator in
\eqref{eq:velocity-operator} unchanged, but avoids repeating part of the
comparator and uncomputation overhead.

\subsection{Pressure circuit}
\label{sec:pressure-circuit}

The pressure circuit starts from the same shared rotated state
\eqref{eq:shared-rotated-forcing}, but follows the longitudinal branch \(c=0\).
Using the pressure relation \eqref{eq:pressure-fourier-solution} in the
longitudinal--transverse basis, pressure recovery requires the scalar block
\(U_\Lambda\) on \(c=0\), followed by the phase \(U_{-i}\) and the inverse QFT.

The postselected pressure operator is
\begin{equation}
	\mathcal U_p
	=
	(\mathrm{QFT}_{2D}^{\dagger}\otimes I_c)
	U_{-i}
	U_\Lambda
	U_R
	(\mathrm{QFT}_{2D}\otimes I_c)
	U_I .
	\label{eq:pressure-operator}
\end{equation}

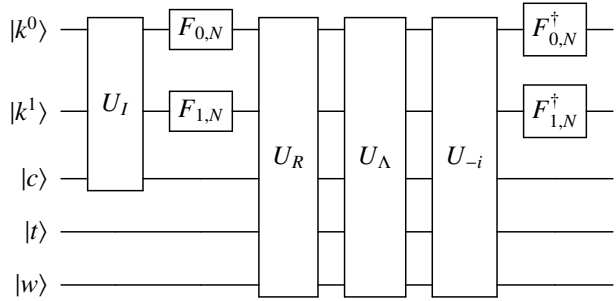
\begin{figure}[h]
	\centering
	\[
	\Qcircuit @C=1.0em @R=1.1em {
		\lstick{\lvert k^0 \rangle} & \multigate{2}{U_I} & \gate{\mathrm{QFT}_{0,N}} & \multigate{4}{U_R} & \multigate{4}{U_\Lambda} & \multigate{4}{U_{-i}} & \gate{\mathrm{QFT}_{0,N}^\dagger} & \qw \\
		\lstick{\lvert k^1 \rangle} & \ghost{U_I}        & \gate{\mathrm{QFT}_{1,N}} & \ghost{U_R}        & \ghost{U_\Lambda}        & \ghost{U_{-i}}        & \gate{\mathrm{QFT}_{1,N}^\dagger} & \qw \\
		\lstick{\lvert c \rangle}   & \ghost{U_I}        & \qw            & \ghost{U_R}        & \ghost{U_\Lambda}        & \ghost{U_{-i}}        & \qw                    & \qw \\
		\lstick{\lvert t \rangle}   & \qw                & \qw            & \ghost{U_R}        & \ghost{U_\Lambda}        & \ghost{U_{-i}}        & \qw                    & \qw \\
		\lstick{\lvert w \rangle}   & \qw                & \qw            & \ghost{U_R}        & \ghost{U_\Lambda}        & \ghost{U_{-i}}        & \qw                    & \qw
	}
	\]
	\caption{Quantum circuit for the pressure field on registers
		\(k^0,k^1,c,t,w\). The pressure uses the shared preprocessing
		\(U_I\), \(\mathrm{QFT}_{2D}\), and \(U_R\), then applies \(U_\Lambda\) and the
		controlled phase \(U_{-i}\) on the longitudinal selected branch
		\(c=0\), \(t=1\), and \(w=0\).}
	\label{fig:pressure-circuit}
\end{figure}

Unlike the velocity circuit, no inverse component rotation is needed: the
pressure output is scalar, and \(c=0\) is part of the postselection predicate
rather than a reconstructed Cartesian degree of freedom
(Figure~\ref{fig:pressure-circuit}).

The phase gate $U_{-i}$ \hl{is defined on the pressure-success branch by}:
\begin{equation}
U_{-i}:\;
|k^0\rangle|k^1\rangle|0\rangle_c|1\rangle_t|0\rangle_w
\mapsto
-i|k^0\rangle|k^1\rangle|0\rangle_c|1\rangle_t|0\rangle_w,
\end{equation}
\hl{and} is controlled on the pressure-success branch
\(c=0,\ t=1,\ w=0\). On that branch it applies the factor \(-i\); on the
orthogonal branches it acts as the identity.

Conditioned on \(c=0\), \(t=1\), and \(w=0\), the final state after
\(\mathrm{QFT}_{2D}^{\dagger}\) is proportional to the physical-space pressure field,
\begin{equation}
	|\psi_p\rangle
	\propto
	\sum_{k^0,k^1}
	p(k^0,k^1)|k^0\rangle|k^1\rangle,
	\label{eq:pressure-success-amplitudes}
\end{equation}
where, to obtain the physical values for pressure, the selected amplitudes should be scaled by 
\(\|\bm f\|/\varepsilon_\Lambda\).

\subsection{Shared preprocessing, separated extraction, and block-encoding alternative}

The velocity and pressure circuits are written separately to expose their
different postselection predicates, but they are not independent algorithms. The
state preparation \(U_I\), the two-dimensional QFT \(\mathrm{QFT}_{2D}\), and the angular
component rotation \(U_R\) are common to both. After this point the extraction
semantics diverge. Velocity uses the transverse branch, the Green block
\(U_\Gamma\), and the inverse component rotation \(U_R^\dagger\). Pressure uses
the longitudinal branch, the scalar block \(U_\Lambda\), and the controlled
phase \(U_{-i}\), with no inverse component rotation.

A single monolithic circuit could contain both extractions, but it would mix
branches with different postselection predicates, different scale factors
\(\varepsilon_\Gamma\) and \(\varepsilon_\Lambda\), and different final
reconstructions. Keeping the pathways separate makes the physical interpretation
and the rescaling transparent: the velocity branch \hl{uses $c=1$ to control $U_\Gamma$} before \(U_R^\dagger\) and
through \(t=1,w=0\) after the scalar block, whereas the pressure branch is
\(c=0,t=1,w=0\). This separation is also the convention
used for circuit-size reporting. In the numerical section, circuits are
transpiled to the universal basis formed by CNOT gates and single-qubit \(U_3\)
rotations, and the corresponding gate counts are used as the implementation-level
complexity measure.

The rotation-based projector is the implementation used in the numerical
examples because it is economical in two dimensions: one angle determines the
change from Cartesian to longitudinal--transverse components. A complementary
linear-combination-of-unitaries (LCU) representation, inspired by the
block-encoding multiscale formulation of Wang et al.~\cite{Wang2026}, is not
used in the benchmarks but is useful for assessing extensions. In three
dimensions, a rotation-based implementation would require a full
mode-dependent orthonormal frame around each nonzero \(\bm k\), with additional
choices and singular directions. The LCU route instead represents the vector
projector directly as a sum of simpler unitary component operations, avoiding an
explicit parameterization of that frame. The LCU form of the velocity projector and the pressure row-contraction alternative is given in~\ref{app:lcu-block-encoding}.

\subsection{Operator-wise complexity estimate}
\label{subsec:operator-complexity}

A theoretical operator-by-operator complexity order can be estimated. Let \(n=\log_2 N\) denote the number of
qubits per spatial direction, and let \(T_R\), \(T_\Gamma\), and \(T_\Lambda\) denote the tile counts used for \(U_R\), \(U_\Gamma\), and \(U_\Lambda\),
respectively. 
Each one-dimensional QFT costs \(O(n^2)\) gates \hl{from the universal basis previously defined}, so the two-dimensional transform contributes \(O(n^2)\) up to constant factors. 

\hl{Following the piecewise polynomial encoding of Liu et al.~\cite{Liu2024}, the bivariate polynomial construction analysed there contains the following number of monomial terms for two $n$-qubit variables and degree $p$:}
\begin{equation}
	\hll{M(n,p)=\binom{n+p}{p}^{2}.}
\end{equation}
\hl{Liu et al. combine this term count with the linear-in-$n$ decomposition cost of multi-controlled rotations to obtain their bivariate polynomial-encoding complexity estimate. Accordingly, a block with $T$ tiles contributes $O(T\,\mathrm{poly}(p,n))$ gates in the present formulation, accounting for tile comparators, controlled polynomial rotations, and uncomputation. For the velocity circuit, this generic block cost applies to $U_R$ with tile count $T_R$ and degree $p_R$, and to $U_\Gamma$ with tile count $T_\Gamma$ and degree $p_\Gamma$. The following conditional gate cost accounts for the spectral operator blocks only; it excludes generic state preparation, postselection, and measurement. The corresponding universal-gate complexity of the velocity circuit is} 

\begin{equation}
	C_{\bm u}
	=
	O(n^2)
	+
	O\!\left(T_R\,\mathrm{poly}(p_R,n)\right)
	+
	O\!\left(T_\Gamma\,\mathrm{poly}(p_\Gamma,n)\right),
\end{equation}
up to constant factors associated with the forward and inverse QFTs and the
repeated use of \(U_R\) and \(U_R^\dagger\). The pressure circuit has the same structure with \(U_\Lambda\) replacing \(U_\Gamma\) and without the inverse component rotation. The transverse selection is absorbed into the
component-conditioned application of \(U_\Gamma\), while the pressure phase
\(U_{-i}\) is a lower-order controlled phase. Thus, for fixed approximation
parameters, the grid-resolution dependence is polylogarithmic in the number of grid points, while tile counts, polynomial degrees, success probabilities, state preparation, and observable-estimation costs remain explicit implementation parameters.

\hl{For comparison, one application of the corresponding classical FFT-based Stokes operator on an $N\times N$ grid costs $O(N^2\log N)$, followed by $O(N^2)$ mode-wise multipliers. For fixed tiling and polynomial degrees, the quantum construction has the polylogarithmic resolution dependence given above; this formal difference is its algorithmic motivation.}

\subsection{Homogenized observables}
\label{sec:homogenized-observable}

The velocity circuit of Section~\ref{sec:velocity-circuit} produces, on its
selected branch and up to the known amplitude-encoding and input-normalization
factors, the projected spectral Stokes response.
For full field reconstruction, the circuit is completed by applying
\(U_R^\dagger\) and the inverse QFT. These final unitary operations are not
always needed. 

In multiscale flow calculations, the fine-scale solve is often
used through homogenized quantities rather than through every resolved degree of
freedom. Two such outputs are the kinetic energy and the Reynolds stress of the
periodic fluctuation. Both are accessible from measurements before the inverse
Fourier transform. In an RVE formulation one may decompose
\(\bm u=\overline{\bm u}+\bm u'\); the periodic fine-scale solve considered
here supplies the mean-zero fluctuation \(\bm u'\), or equivalently the
fluctuating response after subtracting the prescribed coarse-scale mean.

In the physical-grid representation, where the labels \(k^0,k^1\) denote grid
points, the kinetic-energy observable is
\begin{equation}
	K
	=
	\frac{1}{2}\left\langle |\bm u|^2 \right\rangle
	=
	\frac{1}{2N^2}\sum_{k^0,k^1}
	\left(
	u_0(k^0,k^1)^2
	+
	u_1(k^0,k^1)^2
	\right).
	\label{eq:kinetic-energy-definition}
\end{equation}
With the unitary Fourier convention used by \(\mathrm{QFT}_{2D}\), Parseval's identity
allows the same quantity to be evaluated in spectral space. \hl{Here $N^{-2}$ is the physical-grid averaging factor, distinct from the unitary normalization of the QFT coefficients.} Since the selected
transverse response satisfies
\[
\hat u_\perp(\bm k)
=
\Gamma(\bm k)\hat f_\perp(\bm k),
\]
one obtains
\begin{equation}
	K
	=
	\frac{1}{2N^2}
	\sum_{k^0,k^1}
	\Gamma(\bm k)^2
	\left|\hat f_\perp(\bm k)\right|^2 .
	\label{eq:kinetic-energy-parseval}
\end{equation}
Define the shortened observable block
\begin{equation}
	\mathcal A_K
	=
	U_\Gamma U_R
	(\mathrm{QFT}_{2D}\otimes I_c)U_I .
	\label{eq:AK-definition}
\end{equation}
This block stops after \(U_\Gamma\), before applying \(U_R^\dagger\) and the
inverse QFT. On the selected branch,
\begin{equation}
	\frac{\varepsilon_\Gamma}{\|\bm f\|}
	\Gamma(\bm k)\hat f_\perp(\bm k)
	|k^0\rangle|k^1\rangle|1\rangle_c|1\rangle_t|0\cdots0\rangle_w
	\label{eq:kinetic-success-amplitude}
\end{equation}
is the selected contribution associated with each mode label. The kinetic
energy is therefore encoded as a postselection probability of the existing registers.
Let
\begin{equation}
	\Pi_K
	=
	I_{k^0k^1}
	\otimes |1\rangle\langle 1|_c
	\otimes |1\rangle\langle 1|_t
	\otimes |0\cdots0\rangle\langle0\cdots0|_w.
	\label{eq:kinetic-projector}
\end{equation}
Then
\begin{equation}
	P_K
	=
	\left\langle 0\right|
	\mathcal A_K^\dagger \Pi_K \mathcal A_K
	\left|0\right\rangle
	=
	\Pr(c=1,\;t=1,\;w=0\cdots0)
	\label{eq:kinetic-branch-probability}
\end{equation}
is obtained operationally by measuring the component, target, and work
registers after \(\mathcal A_K\) and counting the frequency of the postselection. The mode registers need not be measured. This avoids reconstructing all
Fourier coefficients, but it does not remove the sampling or amplitude-estimation
cost associated with estimating a probability. 
\hl{The grid-level kinetic-energy approximation is recovered from the postselection probability as}
\begin{equation}
	K
	=
	\frac{\|\bm f\|^2}{2N^2\varepsilon_\Gamma^2}\,P_K = \frac{F^2}{2\varepsilon_\Gamma^2}\,P_K.
	\label{eq:kinetic-energy-from-probability}
\end{equation}
On a shot-based backend, estimating \(P_K\) means repeating the circuit and
estimating the frequency of the selected measurement outcome. This is the
standard operational way of estimating a branch probability in the circuit
model. \hl{~\ref{app:postselection-improvements} quantifies this
selected-branch probability, its direct-sampling cost, and two mitigation strategies for the kinetic-energy observable.}
The factor \(\|\bm f\|^2/\varepsilon_\Gamma^2\) restores the physical forcing
normalization and the scalar amplitude-encoding scale. 

The same prepared Stokes state also gives access to the Reynolds stress tensor,
or second moment of the fluctuation,
\begin{equation}
	R_{ij}
	=
	\langle u'_i u'_j\rangle
	=
	\frac{1}{N^d}\sum_{k^0,k^1}
	\hat u'_i(\bm k)^*\hat u'_j(\bm k),
	\qquad i,j=0,\ldots,d-1 .
	\label{eq:fluctuation-tensor-definition}
\end{equation}
where \(i\) and \(j\) are free component indices.
For this tensor the component basis must be Cartesian. The circuit therefore
continues with \(U_R^\dagger\), but still stops before the inverse QFT:
\begin{equation}
	\mathcal A_T
	=
	U_R^\dagger \mathcal A_K .
	\label{eq:AT-definition}
\end{equation}
On the selected branch, it prepares the normalized Cartesian spectral
fluctuation state
\begin{equation}
	|\hat{\bm u}'_{\mathrm{norm}}\rangle
	=
	\frac{1}{\|\bm u'\|}
	\sum_{k^0,k^1}\sum_{i=0}^{d-1}
	\hat u'_i(\bm k)|k^0\rangle|k^1\rangle|i\rangle_c ,
	\qquad
	\|\bm u'\|^2
	=
	\frac{\|\bm f\|^2}{\varepsilon_\Gamma^2}P_K .
	\label{eq:normalized-cartesian-spectral-velocity}
\end{equation}
The tensor components are then component-space expectation values,
\begin{equation}
	R_{ij}
	=
	\frac{\|\bm u'\|^2}{N^d}
	\left\langle\hat{\bm u}'_{\mathrm{norm}}\right|
	I_{\mathcal K}\otimes |i\rangle\langle j|_c
	\left|\hat{\bm u}'_{\mathrm{norm}}\right\rangle .
	\label{eq:tensor-component-measurement}
\end{equation}
The diagonal entries of \(R_{ij}\) can be obtained by measuring the component
qubit in the Cartesian basis, conditioned on the same selected branch
\(t=1,w=0\). Off-diagonal entries require measurements in rotated component
bases, since they correspond to coherences between the two Cartesian component
amplitudes. 

\hl{In the present work, we validate the trace quantity $K=\frac12\operatorname{tr} R$, obtained from the selected-branch probability in Eq.~\eqref{eq:kinetic-energy-from-probability}, together with the Reynolds-stress components $R_{xx}$ and $R_{xy}$. In two dimensions, these quantities determine the symmetric Reynolds-stress tensor.}

\begin{figure}[h]
	\centering
	\begin{subfigure}[t]{0.48\textwidth}
		\centering
		\[
		\Qcircuit @C=0.85em @R=1.0em {
			\lstick{\lvert k^0 \rangle} & \multigate{4}{\mathcal A_K} & \qw     & \qw \\
			\lstick{\lvert k^1 \rangle} & \ghost{\mathcal A_K}        & \qw     & \qw \\
			\lstick{\lvert c \rangle}   & \ghost{\mathcal A_K}        & \meter  & \cw \\
			\lstick{\lvert t \rangle}   & \ghost{\mathcal A_K}        & \meter  & \cw \\
			\lstick{\lvert w \rangle}   & \ghost{\mathcal A_K}        & \meter  & \cw
		}
		\]
		\caption{Kinetic energy.}
		\label{fig:kinetic-observable-circuit}
	\end{subfigure}
	\hfill
	\begin{subfigure}[t]{0.48\textwidth}
		\centering
		\[
		\Qcircuit @C=0.85em @R=1.0em {
			\lstick{\lvert k^0 \rangle} & \multigate{4}{\mathcal A_T} & \qw                 & \qw & \qw \\
			\lstick{\lvert k^1 \rangle} & \ghost{\mathcal A_T}        & \qw                 & \qw & \qw \\
			\lstick{\lvert c \rangle}   & \ghost{\mathcal A_T}        & \gate{B_{ij}}       & \meter & \cw \\
			\lstick{\lvert t \rangle}   & \ghost{\mathcal A_T}        & \meter              & \cw & \cw \\
			\lstick{\lvert w \rangle}   & \ghost{\mathcal A_T}        & \meter              & \cw & \cw
		}
		\]
		\caption{Reynolds stress.}
		\label{fig:tensor-observable-circuit}
	\end{subfigure}
    \caption{\label{fig:observable-circuits} Homogenized observable circuits. The kinetic-energy block
	\(\mathcal A_K\) is measured before \(U_R^\dagger\), so \(c=1\) still denotes the transverse branch. The Reynolds-stress block
	\(\mathcal A_T=U_R^\dagger\mathcal A_K\) returns the component register to the
	Cartesian basis. In panel (b), $B_{ij}$ denotes the component-space measurement block implementing $I_{k^0k^1} \otimes |i\rangle\langle j|_c$ for extraction of the Reynolds-stress component $R_{ij}$. In both cases, the mode registers are not measured.}
	\label{fig:homogenized-observable-circuits}
\end{figure}
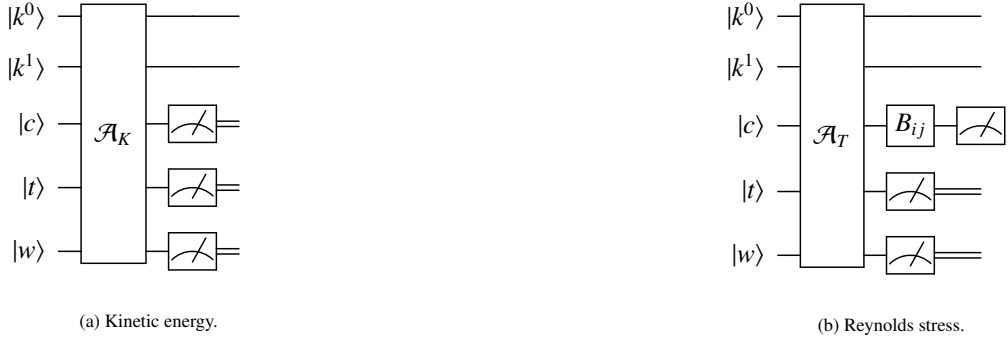

Figure~\ref{fig:observable-circuits} summarizes the two observable-extraction
circuits and highlights where the component register is interpreted as a
transverse branch or as a Cartesian component.

\section{Results and Discussion}~\label{sec:results}

We present three periodic Stokes tests, each targeting a different aspect of
the proposed quantum spectral approach. The Taylor--Green problem is used
as a manufactured verification case to check the consistency of the assembled
velocity and pressure circuits against an analytic solution. The regularized
force-dipole problem is a periodic Stokes benchmark with localized forcing and
broad Fourier support; it is used to assess velocity and pressure convergence
beyond sparse manufactured modes and to distinguish finite-resolution effects
from errors introduced by the polynomial spectral approximation. The
RVE-inspired problem prescribes a periodic fluctuation forcing and focuses on
homogenized outputs, in particular the kinetic-energy observable introduced in
Section~\ref{sec:homogenized-observable}, which can be estimated without
reconstructing the full velocity field. Before these case studies, we first
display the piecewise approximation of the spectral symbols that enter
\(U_R\) and \(U_\Gamma\).

All circuits are implemented in Qiskit and executed on a noiseless state-vector
simulator~\cite{Qiskit2024}. The reported results therefore assess the
algorithmic and approximation properties of the proposed spectral construction,
without modelling hardware noise, finite-shot sampling, or fault-tolerant
overheads. Gate counts are reported after transpilation to the universal gate
set formed by CNOT gates and single-qubit \(U_3\) rotations. They should be read
as implementation-level measurements of the constants and transpilation
overheads associated with the operator-wise scaling discussed in
Section~\ref{subsec:operator-complexity}. \hl{For the polynomial scalar blocks, the velocity and pressure scale factors
were selected separately so that
$\max_{\bm k}|\varepsilon_\Gamma\Gamma(\bm k)|=
\max_{\bm k}|\varepsilon_\Lambda\Lambda(\bm k)|=10^{-2}$ in every reported
configuration. Thus, the standard bound
$|\sin x-x|\le |x|^3/6$ limits the relative small-angle error to less than
$1.7\times10^{-5}$ before the independent polynomial-approximation error is
considered.} Unless stated otherwise, field errors are reported in
the relative discrete (L$^2$) norm,
\begin{equation}
E_{L^2}(q,q_{\mathrm{ref}})
=
\frac{|q-q_{\mathrm{ref}}|_{L^2}}{|q_{\mathrm{ref}}|_{L^2}},
\end{equation}
with vector fields evaluated by concatenating their Cartesian components. For
scalar homogenized observables, the same relative-error convention is used with
the scalar reference value in the denominator.

\subsection{Operator tiling and low-mode resolution}

\begin{figure}[tb]
	\centering
	\includegraphics[width=\textwidth]{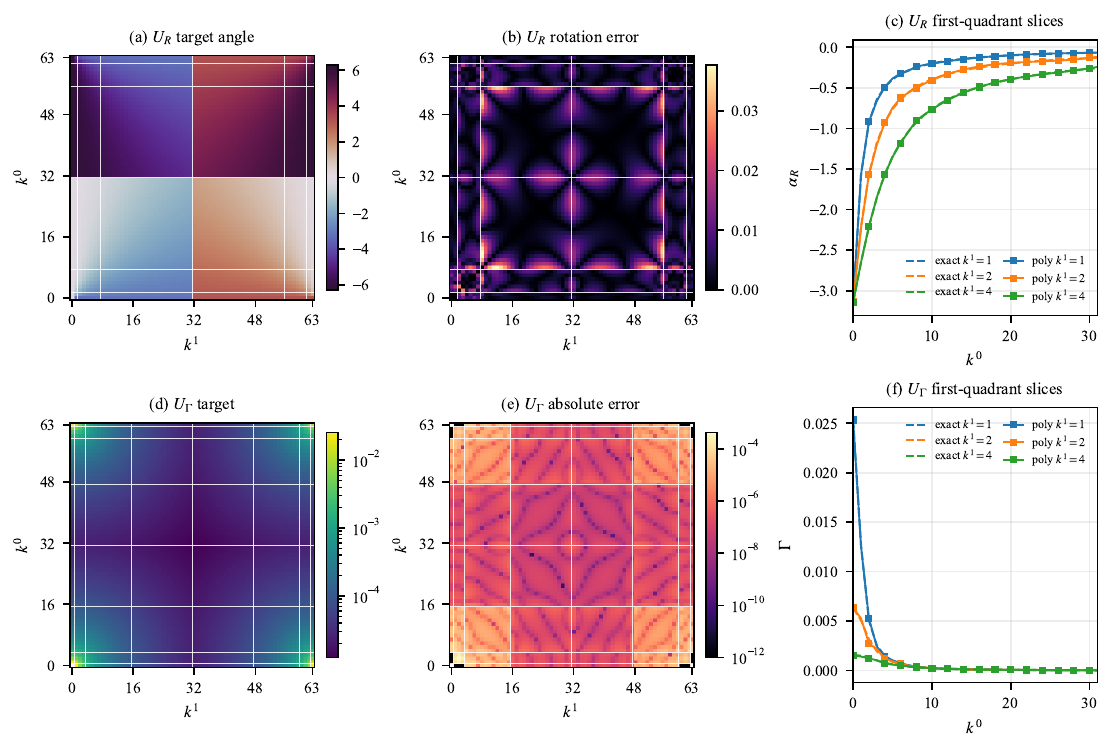}
	\caption{\label{fig:operator-tiling} Piecewise-polynomial operator at \(N=64\). Panels (a)--(c) show the target angle
	\(\alpha_R\), the corresponding approximation error,
	and representative first-quadrant low-mode slices. Panels (d)--(f) show the scalar
	Green factor \(\Gamma\), its absolute approximation error, and the
	corresponding slices. The displayed tiles use the rigid FFT low-mode
	layout used for the Stokes operator. Both maps and polynomial fits
	are shown on raw FFT-index tiles, with degree-three monomial fits in both
	cases.}
\end{figure}

Before considering the three flow cases, it is useful to isolate the spectral
symbols that are approximated by piecewise polynomials. Figure~\ref{fig:operator-tiling}
shows the two maps that dominate the approximation behavior:
the rotation angle \(\alpha_R\) used inside \(U_R\) and the scalar Green factor
\(\Gamma\) loaded by \(U_\Gamma\). In the circuit construction, both the
comparator intervals and the polynomial variables are defined directly in the
raw FFT labels \((k^0,k^1)\). The signed Fourier convention only determines the
symbol value assigned to each raw index. The figure is therefore displayed in
the raw FFT-index plane, where the comparator-based tiling is implemented. The
zero mode is excluded from the Stokes response.

For \(\alpha_R\), the FFT split at \(N/2\) separates the angular branches into
four sectors of the raw index plane. The largest discrete changes occur near
the low-wavenumber indices and across the branch transitions induced by this
split. For \(\Gamma\), there is no angular branch structure, but the sampled
map has its strongest curvature at the raw indices corresponding to the lowest
nonzero signed wavenumbers. The same low-mode logic therefore applies, even
though the scalar map itself is radially symmetric in Fourier space rather than
quadrant-valued.

\hl{The rotation-error diagnostic is also relevant to incompressibility: an error in the approximated angle perturbs the local longitudinal--transverse frame and can generate a residual longitudinal component relative to the exact wavevector after transverse-branch selection. In Section~\ref{sec:results-homogenized}, we quantify these deviations and their influence on observables.}

This tiling is kept rigid as \(N\) changes: the low-mode breakpoints remain
anchored in signed Fourier space and the larger grids add new modes in the
high-wavenumber tail. Thus, increasing \(N\) should not be interpreted as
automatically improving the approximation of already present low modes; that
accuracy is controlled by the local tile and polynomial degree. Increasing
\(N\) mainly expands the represented spectral support. This distinction is used
below when separating finite-resolution truncation from the error introduced by
the polynomial Stokes operator. The dipole benchmark uses a fixed
operator-oriented tiling, while the RVE observable study uses a related
low- and intermediate-wavenumber tiling calibrated to the expected
\(k^{-5/3}\) energy distribution of the populated forcing spectrum.

\subsection{Steady Taylor--Green vortex}

As a first end-to-end check, we use a smooth periodic Taylor--Green
manufactured solution on $\Omega=(0,1)^2$, with $\mu=1$,
\begin{align}
	u_0(x_0,x_1) &= \sin(2\pi x_0)\cos(2\pi x_1), \\
	u_1(x_0,x_1) &= -\cos(2\pi x_0)\sin(2\pi x_1), \\
	p(x_0,x_1) &= \cos(2\pi x_0)\cos(2\pi x_1).
\end{align}
The forcing is generated from ($\bm f=-\mu\Delta\bm u+\nabla p$), sampled on
the periodic grid, normalized into the input quantum state, and compared with
the classical FFT solution on the same grid. This test is used only to verify
the assembled circuit: the Fourier conventions, the pressure phase, the
normalization and reconstruction operators, and the incompressibility projection.
Because the analytic velocity and pressure are known, this case is also useful
for checking pressure recovery and projection semantics in a controlled
calculation.

It is deliberately not used as the main convergence benchmark. In spectral
space, the Taylor--Green forcing is concentrated on a small set of discrete
Fourier modes, effectively a finite collection of modal delta contributions. Once those active modes are represented on the grid, increasing \(N\) or refining the piecewise tiling does not test the approximation of a broader spectral support. The Taylor--Green case is therefore a consistency check for the assembled velocity and pressure circuits, while the force-dipole and RVE-inspired cases below are used to assess broadband spectral approximation.

\subsection{Periodic regularized force-dipole Stokes benchmark}

We next consider a forcing that is native to periodic Stokes flow rather than
manufactured from a single Fourier mode. In microhydrodynamics, a force dipole
is a classical idealization of a small inclusion, an active particle, a local
stirrer, an internal pair of forces, or a localized stress source. Such a
loading injects momentum locally but has no net force. It therefore does not
accelerate the periodic domain as a whole and is compatible with the zero-mode
convention used by the periodic Stokes operator.

The benchmark is defined on the periodic torus by two equal and opposite
regularized forces,
\begin{equation}
\bm f(\bm x)
=
\bm F
\left[
g_{\sigma}^{\mathrm{per}}(\bm x-\bm x_+)
-
g_{\sigma}^{\mathrm{per}}(\bm x-\bm x_-)
\right],
\qquad \bm x=(x_0,x_1),
\end{equation}
where \(\bm F=(1,0)^T\), \(\bm x_+=(0.35,0.50)\),
\(\bm x_-=(0.65,0.50)\), and \(\mu=1\). Here
\(g_{\sigma}^{\mathrm{per}}\) is a grid-normalized periodic Gaussian blob,
implemented using the minimum-image distance on the unit torus. After the two
force components are assembled, the component-wise mean is removed to suppress
the zero Fourier mode up to roundoff. Unlike the Taylor--Green check, this forcing is localized in physical
space and therefore excites a broad range of Fourier modes while remaining
smooth and exactly periodic.

\begin{figure}[tb]
	\centering
	\includegraphics[width=\textwidth]{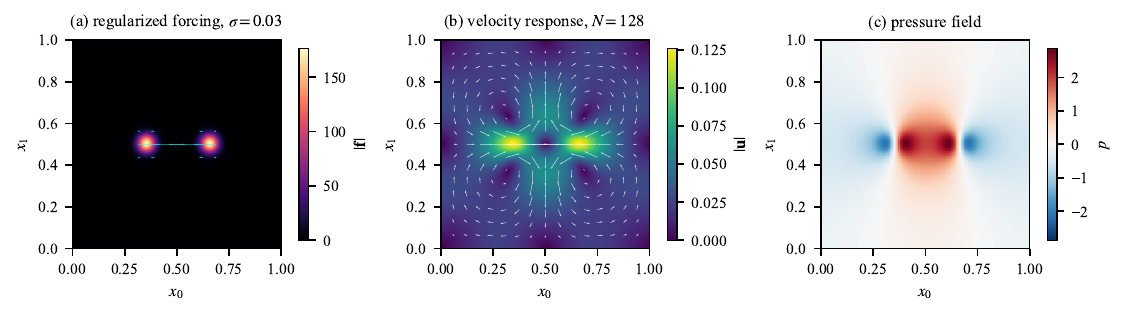}
	\caption{Periodic regularized force-dipole benchmark for \(\sigma=0.03\) and
		\(N=128\). Panel (a) shows the forcing magnitude with direction arrows. Panel
		(b) shows the velocity magnitude and direction. Panel (c) shows the
		corresponding zero-mean pressure field.}
	\label{fig:dipole-fields}
\end{figure}

The parameter \(\sigma\) controls the spectral bandwidth of the benchmark.
Decreasing \(\sigma\) localizes the applied force in physical space and spreads
its Fourier content over a wider range of modes. Thus, the three values of
\(\sigma\) used below define a controlled sequence of spectral difficulty:
large \(\sigma\) gives a smoother, low-bandwidth forcing, whereas small
\(\sigma\) activates higher modes and makes the fixed polynomial tiling more
demanding. This makes the dipole benchmark a useful intermediate test between
the sparse Taylor--Green verification and the broadband RVE-inspired
fluctuation considered in the next section. Figure~\ref{fig:dipole-fields}
shows the forcing, velocity response, and pressure field for a representative
case. The pressure is included to illustrate the longitudinal response induced
by the localized forcing; the error study below reports both velocity and
pressure.

The rotation-angle map \(\alpha_R\) used by \(U_R\) and the scalar Green factor
\(\Gamma\) loaded by \(U_\Gamma\) are implemented with degree-three piecewise polynomial encodings on the same
nested low-wavenumber tiling of the signed Fourier plane. The tiling is
refined near low modes and coarsened toward the Nyquist range. The same tiling
and polynomial degree are used for all values of \(\sigma\), so that decreasing
\(\sigma\) tests the fixed approximation under increasing spectral bandwidth.
The reference solution for the reported errors is the classical FFT solution
computed on a \(512^2\) grid and sampled down to the working grid.
The piecewise layout used here was selected by a trial-and-error calibration
against the target spectral functions. A systematic design of this layout is
left outside the present study. In principle, the tile partition and polynomial
degrees could be optimized by a multiobjective criterion that maximizes the
quality of the operator approximation while minimizing the resulting gate
count.

\begin{figure}[bt]
	\centering
	\includegraphics[width=\textwidth]{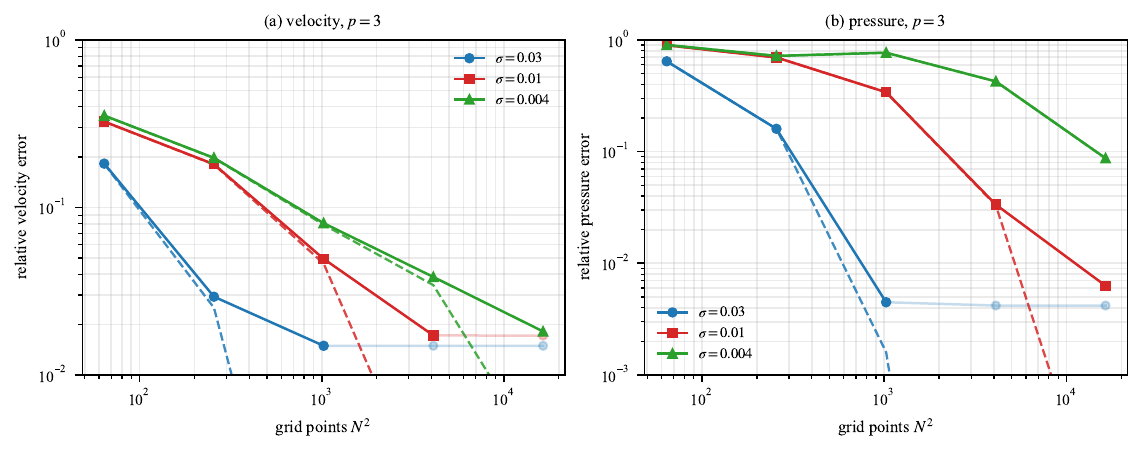}
	\caption{\label{fig:dipole-convergence} Regularized dipole benchmark against the \(512^2\) FFT reference. The
		sweep uses \(\sigma\in\{0.03,0.01,0.004\}\) and
		\(n\in\{3,4,5,6,7\}\). Dashed curves use
		\(q=q_N^{\mathrm{FFT}}\) and show the convergence of the classical implementation of the spectral method. Solid curves use \(q=q_N^{\mathrm{poly}}\) and
		show the error obtained with the configured piecewise polynomial Stokes
		operator.}
\end{figure}

Figure~\ref{fig:dipole-convergence} reports the relative \(L^2\) velocity and
pressure errors against a \(512^2\) FFT reference, sampled down to the
computational grid. Two comparisons are shown. The dashed curves use
\(q=q_N^{\mathrm{FFT}}\) and \(q_{\mathrm{ref}}=q_{512}^{\mathrm{FFT}}\);
they therefore measure the finite-resolution truncation error of the exact
discrete spectral Stokes solver, obtained with a classical implementation of the spectral method. The solid curves use
\(q=q_N^{\mathrm{poly}}\) and the same reference
\(q_{\mathrm{ref}}=q_{512}^{\mathrm{FFT}}\); they show the error obtained with
the piecewise polynomial Stokes operator implemented in the quantum circuit.
The gap between the solid and dashed curves is therefore a diagnostic of the
additional error introduced by the polynomial approximation of the spectral
Stokes operator, on top of the finite-resolution truncation error.

The separation between these two effects is most visible for the smaller values
of \(\sigma\). As \(\sigma\) decreases, more high-frequency modes carry
non-negligible forcing amplitude, while the tiling and polynomial degree are
kept fixed. The exact FFT baseline then measures how much of the localized
forcing is resolved by the grid, whereas any plateau or slower decay of the
solid curves indicates the limitation of the fixed polynomial representation of
the rotation-angle and Green-factor maps. This is the intended role of the
dipole benchmark: it does not test a single active Fourier mode, but probes how
the assembled Stokes block behaves as the forcing bandwidth is increased in a
controlled way.

\subsection{RVE-inspired periodic fluctuation benchmark}

The main numerical example is designed to connect the present Stokes solver
with the RVE setting that motivates quantum spectral methods in computational mechanics. In multiscale
fluid formulations, the fine-scale contribution is commonly represented as a
periodic fluctuation with the mean component removed, and the coarse solver uses
averaged outputs such as fluctuation stresses or kinetic energy. The present
Stokes problem is not a replacement for a Navier--Stokes RVE solve, nor a DNS
time step. It is a linear, periodic, divergence-free spectral block that
isolates one essential operation needed by future quantum RVE workflows:
resolving fine-scale incompressible fluctuations and extracting averaged
quantities.

\subsubsection{Forcing and solution samples}

We prescribe a periodic fluctuation field \(\bm u'\) with zero mean and
divergence-free Fourier content,
\begin{equation}
	\nabla\cdot \bm u'=0,
	\qquad
	\langle \bm u'\rangle=0,
	\label{eq:rve-fluctuation-constraints}
\end{equation}
and assign its modal energy over a finite spectral band according to
\begin{equation}
	E_u(k)\sim k^{-5/3},
	\qquad
	k_{\min}\leq |k|\leq k_{\max}.
	\label{eq:rve-energy-spectrum}
\end{equation}
To generate one realization of this fluctuation, each nonzero Fourier mode in
the prescribed band is assigned a random phase and an amplitude consistent with
the target energy spectrum. The velocity coefficient is then projected onto the
transverse direction so that
\(\bm k\cdot\hat{\bm u'}(\bm k)=0\), and conjugate symmetry is imposed to
obtain a real-valued physical field. Modes outside the band, together with the
zero mode, are set to zero. The random seed is fixed throughout the numerical
study, so that changing \(N\) refines the representation of the same underlying
fluctuation rather than generating a new realization.

The Stokes right-hand side is then assembled so that \(\bm u'\) is the exact
periodic response of the continuous spectral problem,
\begin{equation}
	\hat{\bm f}(\bm k)
	=
	\mu |\bm k|^2\hat{\bm u'}(\bm k),
	\qquad
	\hat p(\bm k)=0 .
	\label{eq:rve-forcing-construction}
\end{equation}

Thus the forcing should be interpreted as a periodic, zero-mean source whose
Stokes response is a prescribed fluctuation field. In the reported runs
\(k_{\min}=2\), \(k_{\max}=127\), \(\mu=1\), and the reference solution is
computed classically on a \(512^2\) grid. The affine part associated with a
coarse RVE loading can be added in postprocessing and is not encoded in the
periodic fluctuation state.

This synthetic construction deliberately gives access to a reference
fluctuation and its spectrum, but the quantum circuit is not given
\(\bm u'\).
It receives only the forcing state and must recover the incompressible response through the same \(U_R\) and \(U_\Gamma\) blocks used in the general Stokes solver. The
purpose is therefore not to generate turbulence with a linear operator, nor to
claim that the spectrum was unknown in this benchmark. The purpose is to test a
controlled Stokes block under a physically motivated broadband spectrum. In a
multiscale flow calculation, analogous right-hand sides would be assembled from
temporal, convective, closure, or coarse-scale coupling terms, and would not be
available as the spectrum of a preselected solution.

\begin{figure}[bt]
\centering
\includegraphics[width=\textwidth]{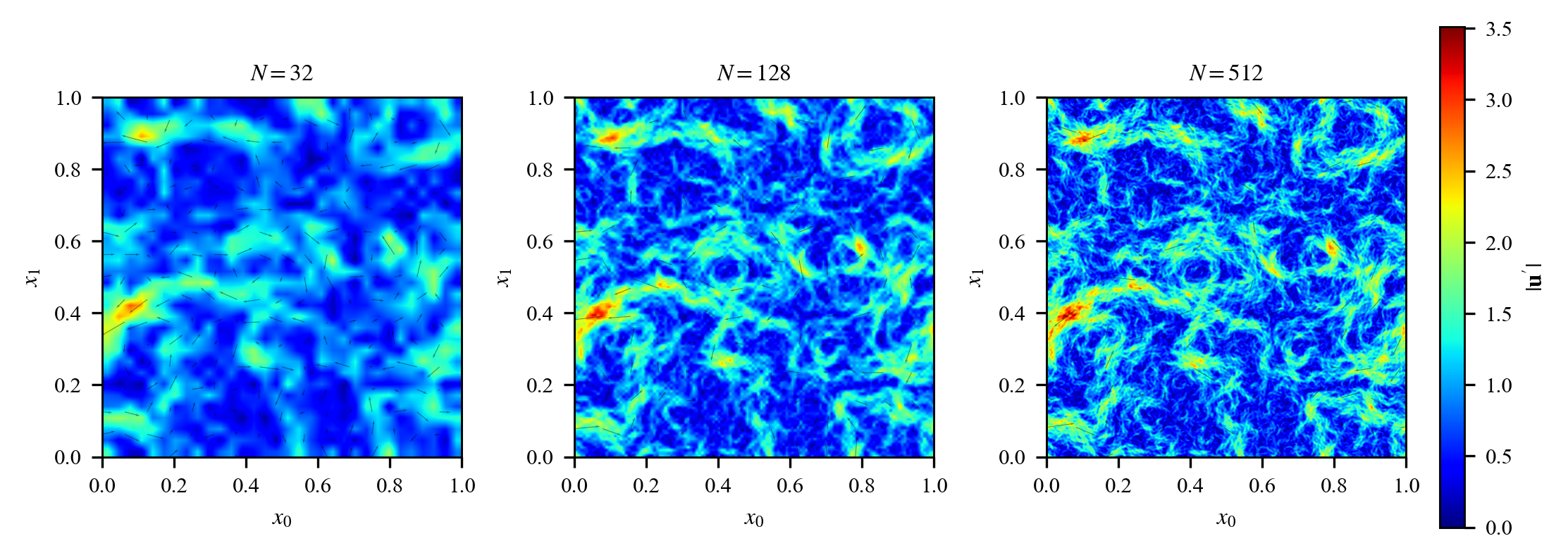}
\caption{RVE-inspired periodic fluctuation fields computed for \(N=32\), \(128\), and \(512\). Colors show \(|\bm u'|\). Increasing \(N\)
incorporates finer Fourier content of the same prescribed periodic fluctuation.}
\label{fig:rve-velocity-fields}
\end{figure}

Figure~\ref{fig:rve-velocity-fields} shows the reconstructed fluctuation field
at increasing grid resolutions. The same underlying periodic fluctuation is
being sampled with increasing spectral bandwidth: moving from \(N=32\) to
\(N=128\) and \(N=512\) progressively incorporates finer-scale modal content.
The refinement is therefore interpreted as resolving a larger portion of the
same fine-scale fluctuation spectrum, rather than as changing the RVE loading.

The RVE runs use the same spectral solver architecture as the previous case
study, but specialize the nested low-wavenumber tiling described in
Figure~\ref{fig:operator-tiling}. For \(\Gamma\), we retain a degree-three
rigid low-mode approximation, since the scalar factor \(1/(\mu|\bm k|^2)\) is
most sensitive at low wavenumbers. For \(\alpha_R\), the relevant approximation
error is not a uniform error over the Fourier plane, but the
error after the resulting rotation acts on the populated input spectrum
\(\widehat{\bm f}(\bm k)\).
We therefore use a low- and intermediate-wavenumber partition calibrated to the
prescribed \(k^{-5/3}\) energy slope of the RVE forcing. This differs from the
dipole calibration and is intended to reflect the type of energy-weighted
accuracy needed in future applications to turbulent-flow fluctuations. Once
chosen, the tiling is kept fixed in the reported sweep; only the polynomial
degree \(p_R\) of the \(\alpha_R\) approximation is varied.

\subsubsection{Homogenized observables}\label{sec:results-homogenized}

A quantity of interest is the homogenized kinetic-energy-like observable $K$.
This choice mirrors the role of averaged RVE outputs in multiscale workflows:
the solver need not reproduce every high-frequency Fourier coefficient with
the same absolute accuracy if the target observable is energy weighted. For a
Kolmogorov-like spectrum, high-frequency approximation errors contribute in
proportion to the smaller amount of energy carried by those modes. This makes
low- and intermediate-frequency accuracy more relevant than uniform accuracy
over the entire Fourier grid.

\begin{figure}[bt!]
\centering
\includegraphics[width=\textwidth]{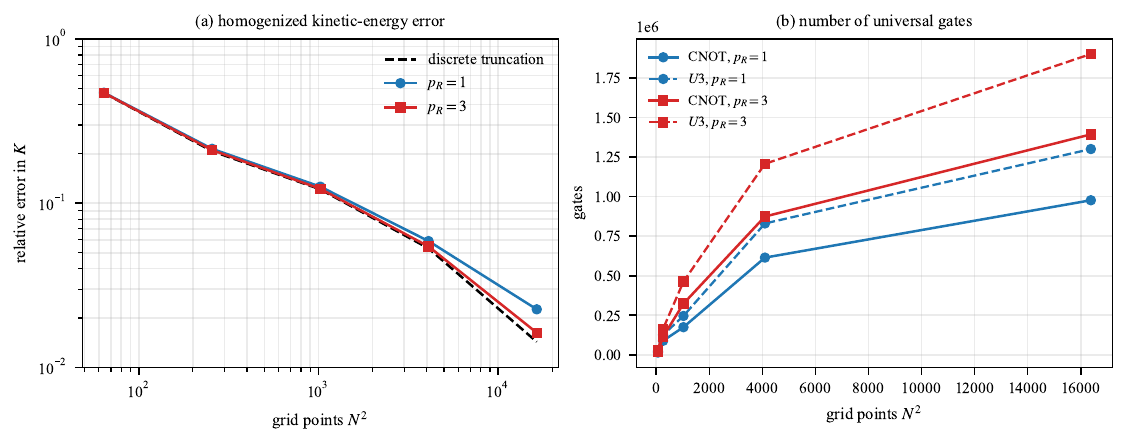}
\caption{RVE-inspired periodic fluctuation benchmark. Panel (a)
reports the relative error of the homogenized kinetic energy
\(K=\frac12\langle |\bm u'|^2\rangle\) against the \(512^2\) classical reference
for rotation-polynomial degrees \(p_R=1\) and \(p_R=3\); the dashed curve is the
discrete truncation error of the exact classical \(N^2\)-grid operator. Panel
(b) reports the number of universal CNOT and \(U_3\) gates after transpilation
for the same fixed nested low-wavenumber tiling.}
\label{fig:rve-results}
\end{figure}

\hl{Figure~\ref{fig:rve-results}(a) shows that the error in the homogenized observable decreases toward the reference value as the grid is refined. The dashed curve is the discrete truncation error, which is computed by applying the exact classical spectral Stokes operator on an $N^2$ grid and comparing the resulting homogenized kinetic energy with the $512^2$ classical reference. It therefore measures finite-resolution truncation of the prescribed spectrum. Then, the comparison with the discrete truncation curve separates the effect of representing more of the prescribed spectrum from the additional error introduced by the polynomial approximations. The results confirm that the approximation error given by $U_R$ can be reduced by increasing $p_R$.}

\begin{figure}[tb!]
	\centering
	\includegraphics[width=\textwidth]{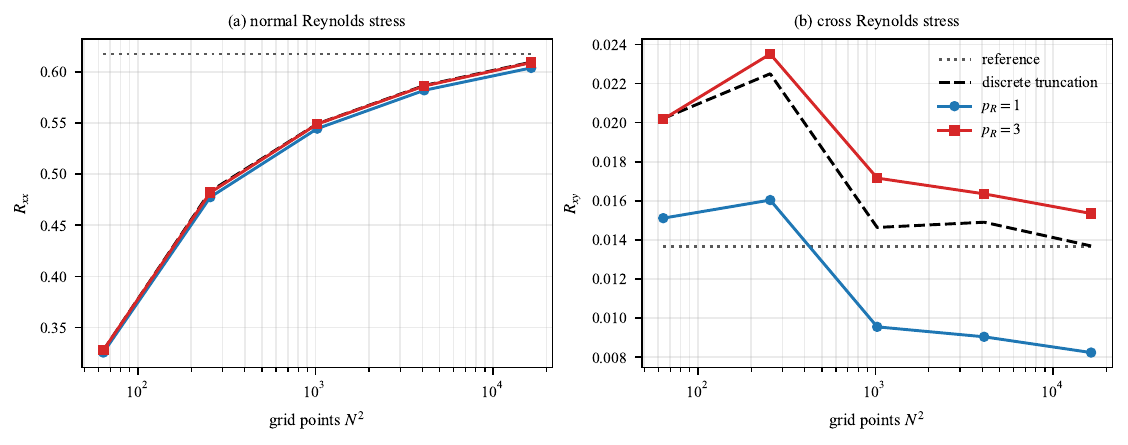}
	\caption{\hl{Reynolds-stress observables for the RVE-inspired periodic fluctuation benchmark. Panel (a) reports the normal component $R_{xx}$, and panel (b) the cross component $R_{xy}$, for the same $N=8,16,32,64,128$ configurations used in Figure~\ref{fig:rve-results}. The dotted line is the fixed $512^2$ FFT reference; the dashed curve is the exact discrete $N^2$-grid operator; solid curves use rotation-polynomial degrees $p_R=1$ and $p_R=3$.}}
\label{fig:rve-reynolds-stress}
\end{figure}

\hl{To extend the RVE validation beyond the trace quantity $K$, we evaluate the Reynolds-stress components $R_{xx}$ and $R_{xy}$ from the same conditions used in Figure~\ref{fig:rve-results}. Figure~\ref{fig:rve-reynolds-stress} compares their values with the reference and with the obtained with the exact discrete $N^2$-grid operator. The normal component $R_{xx}$ approaches the reference as the grid is refined, and the $p_R=3$ rotation approximation is consistently closer to the discrete result than the  $p_R=1$ approximation.}

\hl{The cross component $R_{xy}$ illustrates that convergence of homogenized observables need not be monotone. Unlike energy and normal stresses, which collect non-negative modal contributions, $R_{xy}$ is a signed sum of modal cross contributions; newly resolved modes can therefore enhance or cancel the accumulated value. The discrete truncation nevertheless tends toward the reference as $N$ increases. For $N\geq 32$, increasing $p_R$ from one to three improves the agreement of the polynomial rotation approximation, while the discrete curve isolates the convergence obtained by increasing the number of represented modes.}

\hl{As a postprocessing diagnostic of the reconstructed final Cartesian velocity, we evaluate its divergence relative to the velocity-gradient norm and its relative spectral longitudinal norm. Both quantities would vanish for an exact $U_R$, since the scalar $U_\Gamma$ block preserves the transverse channel. For the fixed tiling, they are not expected to converge monotonically with $N$, because increasing the resolution also expands the high-wavenumber region represented by the outer polynomial tiles. At fixed $N$, however, increasing the rotation degree from $p_R=1$ to $p_R=3$ reduces the normalized divergence residual from $2.0$--$6.7\%$ to at most $1.4\%$ and the longitudinal projection defect from $3.9$--$5.8\%$ to at most $1.1\%$. This reduction accompanies the improved polynomial--discrete agreement of $K$, $R_{xx}$, and $R_{xy}$ in Figures~\ref{fig:rve-results} and~\ref{fig:rve-reynolds-stress}; because the tiling and $U_\Gamma$ approximation are unchanged between the two $p_R$ cases, their differences isolate the benefit of improving the $U_R$ angle fit.}

\subsubsection{Gate-count scaling}

\begin{table}[tb]
\centering
\begin{minipage}[t]{0.47\textwidth}
\centering
\captionof{table}{\hl{Operator-wise gate decomposition for the velocity circuit at $N=128$, transpiled to CNOT and $U_3$ gates.}}
\label{tab:rve-stage-gates}
\begin{tabular*}{\linewidth}{@{\extracolsep{\fill}}lrrr@{}}
\toprule
Stage & CNOT & \(U_3\) & Total \\
\midrule
\(U_I\) & 32\,752 & 32\,767 & 65\,519 \\
QFT & 878 & 1\,384 & 2\,262 \\
\(U_R\) (\(p_R=3\)) & 442\,176 & 647\,580 & 1\,089\,756 \\
\shortstack{\(U_\Gamma\)} & 916\,556 & 1\,218\,077 & 2\,134\,633 \\
\bottomrule
\end{tabular*}
\end{minipage}\hfill
\begin{minipage}[t]{0.49\textwidth}
\centering
\captionof{table}{\hl{Relative transpiled universal-gate count of the isolated $U_R$ block, where $\#U_R(p_R)$ denotes its CNOT+$U_3$ count. Entries are normalized by the $p_R=1$ value at fixed resolution, using the RVE configuration and prescribed tiling of Figure~\ref{fig:rve-results}.}}
\label{tab:rve-degree-sweep}
\begin{tabular*}{\linewidth}{@{\extracolsep{\fill}}rrrrr@{}}
\toprule
\shortstack{Grid points (\(N^2\))} & \shortstack{\(\#U_R(1)\)} & \shortstack{\(\frac{\#U_R(2)}{\#U_R(1)}\)} & \shortstack{\(\frac{\#U_R(3)}{\#U_R(1)}\)} & \shortstack{\(\frac{\#U_R(4)}{\#U_R(1)}\)} \\
\midrule
256 & 15\,904 & 2.22 & 5.06 & 8.11 \\
1\,024 & 49\,452 & 2.69 & 8.38 & 19.09 \\
4\,096 & 61\,572 & 2.98 & 11.33 & 31.82 \\
16\,384 & 73\,692 & 3.28 & 14.79 & 49.62 \\
\bottomrule
\end{tabular*}
\end{minipage}
\end{table}

\hl{Figure~\ref{fig:rve-results}(b) reports the CNOT and $U_3$ counts of the assembled velocity circuit for fixed $p_R=1$ and $p_R=3$. Across the RVE sweep and for the prescribed tiling, the total gate counts are empirically consistent with the polylogarithmic dependence on the $N^2$ grid points predicted by the conditional operator-wise complexity estimate of Section~\ref{subsec:operator-complexity}.}

\hl{To identify how this total cost is distributed among the assembled blocks, Table~\ref{tab:rve-stage-gates} gives the operator-wise decomposition at the largest reported resolution.} The dominant block is $U_\Gamma$, because the RVE Green-factor approximation uses the fixed $64$-tile layout with a degree-three polynomial on each tile.

\hl{To examine the degree dependence associated with the term-count estimate of Section~\ref{subsec:operator-complexity}, Table~\ref{tab:rve-degree-sweep} isolates the $U_R$ contribution from the total-circuit trend of Figure~\ref{fig:rve-results}(b). Its first numerical column reports the $p_R=1$ resolution dependence, whereas each row reports the variation with $p_R$ at fixed resolution. At fixed $p_R$, the isolated $p_R=1$ values in Table~\ref{tab:rve-degree-sweep} confirm the polylogarithmic grid-point dependence predicted by polynomial encoding. 
Conversely, at each fixed resolution, Table~\ref{tab:rve-degree-sweep} shows a monotone and superlinear increase of the $U_R$ gate count over the investigated range $p_R=1,\ldots,4$. This behaviour is consistent with the polynomial degree dependence anticipated theoretically.}

\subsection{Discussion}

The three numerical experiments separate algebraic verification from spectral
approximation. The Taylor--Green manufactured solution verifies the assembled
velocity and pressure pathways, including Fourier conventions, normalization
factors, pressure phase, and the interpretation of the transverse and
longitudinal branches. Because this forcing is supported on only a few Fourier
modes, it should not be interpreted as evidence of broadband convergence of the
polynomial spectral blocks.

The regularized force-dipole benchmark provides the stricter full-field test.
As \(\sigma\) decreases, the physical forcing becomes more localized and its
Fourier support broadens. The comparison with the exact discrete spectral
baseline separates finite-resolution truncation from the additional error of
the polynomial rotation and Green-operator blocks. When the dashed baseline has already decreased but the polynomial curves plateau, the remaining discrepancy is attributable primarily to the fixed-degree, fixed-tiling approximation of the mode-dependent Stokes symbol.

The RVE-inspired benchmark tests a different use case: extracting an averaged
quantity from a fine-scale incompressible fluctuation. The kinetic-energy
observable approaches the refined-grid reference because the dominant
contribution comes from low and intermediate modes, where the energy-weighted
tiling resolves the angular and inverse-Laplacian operators most accurately. This
supports the intended use of the method as a quantum spectral block for
homogenized observables, rather than as a full-field reconstruction tool.

\hl{This is also the relevant interface for a future multiscale use: macroscopic variables can parameterize structured microscopic conditions, and the RVE block can return averaged quantities without requiring pointwise reconstruction of the fine-scale field. State preparation in such an application need not consist of independently loading an arbitrary pointwise field for every RVE. Instead, the microscopic inputs may belong to structured families specified by reduced macroscopic data; when these families admit efficient functional representations, polynomial amplitude encodings provide a route to their coherent preparation.}

\hl{Beyond the RVE setting, the same separation is relevant to time-dependent incompressible flow. At the operator level, a periodic fractional-step method contains a linear substep in which the longitudinal--transverse decomposition can be reused: pressure-correction and viscous updates are represented by scalar Fourier multipliers set by the chosen time discretization (e.g., $(1+\nu\Delta t|\bm k|^2)^{-1}$ for a semi-implicit viscous step).}

The reported gate counts should be interpreted as implementation-level
diagnostics for the transpiled spectral blocks. \hl{~\ref{app:postselection-improvements} quantifies the selected kinetic-energy branch, its ideal direct-sampling cost, and the trade-offs of arcsine-angle encoding and amplitude amplification. These selected-branch estimates complement the transpiled gate counts, while a complete end-to-end resource assessment additionally requires state preparation, reflection synthesis, complete-field readout where needed, hardware noise, and fault-tolerant compilation.}

The main numerical limitation is the construction of
accurate and resource-efficient polynomial encodings for the angular projector
and scalar Green operators. Designing adaptive or energy-weighted tilings is the
most direct route to improving broadband and RVE-oriented performance.

\section{Conclusions}~\label{sec:conclusions}

We have presented a spectral quantum formulation and circuit realization for the steady periodic two-dimensional incompressible Stokes system. The central point is that the periodic Stokes operator admits a compact spectral factorization into a scalar Green operator and a Helmholtz projection. In two dimensions this projection can be implemented through a single mode-dependent component rotation, followed by component-conditioned scalar amplitude blocks for the transverse velocity response and the longitudinal pressure response. This gives a vector-valued incompressible extension of previous quantum spectral methods in computational mechanics while keeping the role of each unitary, postselection predicate, and rescaling factor explicit.

The numerical examples confirm this construction in three complementary
settings: a Taylor--Green manufactured solution for end-to-end circuit
verification, a regularized force-dipole problem for broadband field recovery,
and an RVE-inspired Kolmogorov-like fluctuation for observable extraction. The
last case shows how the circuit can recover the homogenized kinetic energy
\(K=\frac12\langle |\bm u'|^2\rangle\) without reconstructing the full velocity
field, which is the relevant access pattern in many multiscale calculations.

\hl{The proposed solver is a Stokes spectral block for quantum computational mechanics, providing the incompressible linear response and observable-extraction operations needed in a quantum--classical RVE architecture. In such a setting, coarse-scale variables can parameterize microscopic conditions and the RVE responses can be used through homogenized quantities, preserving the association with their corresponding macroscopic data.}

Several limitations remain explicit. The reported gate counts are
implementation-level diagnostics for noiseless state-vector simulations, not
end-to-end runtime estimates. \hl{\ref{app:postselection-improvements} quantifies the selected kinetic-energy branch through ideal direct sampling, arcsine-angle encoding, and amplitude amplification. A complete resource assessment nevertheless still requires state preparation, controlled-reflection synthesis, complete-field readout where required, hardware noise, and fault-tolerant compilation overheads.} The main approximation challenge is not the Fourier transform itself,
but the accurate and resource-efficient polynomial encoding of the angular
projector and scalar Green operators.

The extension path toward full incompressible Navier--Stokes is also clear at
the operator level. A projection or fractional-step method would reuse the
Stokes block as the linear viscous--pressure solve inside each time step, while
the nonlinear convective term could be assembled by a pseudospectral product,
convolutional block encoding, or a hybrid quantum--classical update. Alternative
linearizations, including Carleman embeddings, offer another route for treating
nonlinearity within a quantum linear-system or block-encoding framework. In all
cases, the central requirements are to preserve incompressibility, control the
normalization and postselection costs of non-unitary updates, and extract
coarse-scale observables without unnecessary full-field readout. Generalization
to three-dimensional, non-periodic, and geometry-dependent flows will likely
favor LCU or block-encoding formulations of the vector projector, together with
quantum sine/cosine transforms or domain-extension strategies for boundary
conditions.

\section{Acknowledgments}~\label{sec:ack}
J. Gimenez acknowledges support from MCIU/AEI/10.13039/501100011033 and
FSE+ through project RYC2023-044023-I. 

M. Ortiz gratefully acknowledges the financial support of the Centre Internacional de M\`etodes Num\`erics a l’Enginyeria (CIMNE) of the Universitat Polit\`ecnica de Catalunya (UPC), Spain, through the UNESCO Chair in Numerical Methods in Engineering. 

\appendix
\setcounter{figure}{0}
\setcounter{table}{0}

\section{LCU and block-encoding alternative}
\label{app:lcu-block-encoding}

The rotation-based circuit used in the numerical examples diagonalizes the
two-dimensional Helmholtz projector by a mode-dependent change of component
basis. A complementary construction is possible by expanding the projector
directly as a linear combination of fixed component-space unitaries. It is not used
in the benchmarks of this paper, but it is useful because it avoids constructing
a full mode-dependent orthonormal frame, which becomes more cumbersome in three
dimensions.

\subsection{Pauli decomposition of the Stokes projector}

For every nonzero mode, using \(r_0=r(k^0)\), \(r_1=r(k^1)\), and
\(|\bm k|^2=r_0^2+r_1^2\), the transverse projector is
\begin{equation}
P_\perp(\bm k)
=
I-\frac{1}{|\bm k|^2}
\begin{pmatrix}
r_0^2 & r_0r_1\\
r_0r_1 & r_1^2
\end{pmatrix}.
\end{equation}
In the component basis \(\{|0\rangle_c,|1\rangle_c\}\), this projector can be
written as
\begin{equation}
P_\perp(\bm k)
=
\alpha_0(\bm k)I
+\alpha_X(\bm k)X
+\alpha_Z(\bm k)Z,
\end{equation}
where
\begin{equation}
\alpha_0(\bm k)=\frac12,\qquad
\alpha_X(\bm k)=-\frac{r_0r_1}{|\bm k|^2},\qquad
\alpha_Z(\bm k)=\frac{r_1^2-r_0^2}{2|\bm k|^2}.
\end{equation}
The zero mode is assigned zero response, as in the rotation-based
implementation.

\subsection{Velocity operator as an LCU}

The spectral velocity operator can be written as
\begin{equation}
\Gamma(\bm k)P_\perp(\bm k)
=
\sum_{\ell\in\{0,X,Z\}}
\beta_\ell(\bm k)V_\ell,
\qquad
V_\ell\in\{I,X,Z\},
\end{equation}
with
\begin{equation}
\beta_0=\Gamma\alpha_0,\qquad
\beta_X=\Gamma\alpha_X,\qquad
\beta_Z=\Gamma\alpha_Z.
\end{equation}
Each coefficient \(\beta_\ell(\bm k)\) is a scalar spectral function and can be
encoded with the same piecewise polynomial amplitude-loading machinery used for
\(U_\Gamma\). The difference is that the component-space action is now a
select operation over the fixed unitaries \(I\), \(X\), and \(Z\), rather than a
single angle-dependent \(R_y\) rotation.

Let \(a\) be an LCU ancilla that stores the three labels. A standard
prepare--select--unprepare construction is
\begin{equation}
U_{\mathrm{prep}}|0\rangle_a
=
\sum_{\ell\in\{0,X,Z\}}\sqrt{\omega_\ell}\,|\ell\rangle_a,
\qquad
\omega_\ell>0,\qquad
\sum_\ell\omega_\ell=1,
\end{equation}
\begin{equation}
U_{\mathrm{select}}
=
\sum_{\ell\in\{0,X,Z\}}
|\ell\rangle\langle\ell|_a
\otimes U_{\beta_\ell}
\otimes V_\ell ,
\end{equation}
and
\begin{equation}
U_{\mathrm{LCU}}
=
U_{\mathrm{prep}}^\dagger
U_{\mathrm{select}}
U_{\mathrm{prep}} .
\end{equation}
Postselecting the LCU ancilla back to \(|0\rangle_a\), together with the
successful amplitude-encoding branch of \(U_{\beta_\ell}\), block-encodes an
operator proportional to
\(\Gamma(\bm k)P_\perp(\bm k)\). Signed coefficients can be represented by
separating sign and magnitude in \(U_{\beta_\ell}\), or equivalently by adding
controlled phase factors in the selected branch.

\subsection{Pressure and longitudinal LCU form}

The pressure circuit in the body uses the rotation \(U_R\) to expose
\(\hat f_\parallel(\bm k)\), after which only the scalar operator
\(\Lambda(\bm k)=1/|\bm k|\) and the fixed phase \(-i\) are needed. In an LCU
formulation one may instead write the longitudinal projector as
\begin{equation}
P_\parallel(\bm k)
=
I-P_\perp(\bm k)
=
\frac12 I-\alpha_X(\bm k)X-\alpha_Z(\bm k)Z .
\end{equation}
The longitudinal block entering pressure extraction is then represented by the
same selected component unitaries with scalar maps
\begin{equation}
\Lambda(\bm k)
\left\{
\frac12,\,
-\alpha_X(\bm k),\,
-\alpha_Z(\bm k)
\right\}.
\end{equation}
After the LCU block has isolated the longitudinal component, the successful
pressure branch receives the fixed phase \(U_{-i}=P(-\pi/2)\), and the zero mode
is suppressed to impose the zero-mean pressure gauge. This pressure LCU is
algebraically parallel to the velocity LCU; the difference is the longitudinal
coefficient set and the scalar phase convention.

\subsection{Implementation notes}

The full LCU circuit would use the same forcing preparation \(U_I\), the same
two-dimensional QFT \(\mathrm{QFT}_{2D}\), and the same polynomial approximation backend
as the rotation-based solver. Its additional registers are the LCU label
ancilla and the success registers required by the scalar coefficient encoders.
The extraction semantics are also similar: the useful block is obtained by
postselecting the LCU label ancilla after unpreparation, the scalar
amplitude-encoding success flag, and any work registers returned to zero. A full description of the implementation of this LCU/block-encoding approach is given by Wang et al.~\cite{Wang2026}.

In two dimensions the rotation-based circuit is more compact because the
projector is represented by one angle. The LCU construction becomes more
attractive in three dimensions. There, the transverse projector remains the
simple matrix \(I-\bm k\bm k^T/|\bm k|^2\), but a rotation-based approach would
need a mode-dependent orthonormal frame with additional gauge choices and
discontinuities. Expanding the \(3\times3\) projector in a fixed component-space
unitary basis moves that complexity into scalar coefficient maps, which can be
handled by the same piecewise polynomial machinery already used for
\(\Gamma\), \(\Lambda\), and the two-dimensional rotation angle.

\section{Success probabilities, sampling, and practical mitigation strategies}
\label{app:postselection-improvements}

\hl{The practical usability of a circuit whose output is recovered from a
	selected measurement event is governed in part by the probability of that
	event. This appendix relates that probability to the repetitions required to
	estimate it at a prescribed tolerance, examines its behavior under spectral
	refinement, and records two practical mitigation strategies. For independent
	shots, the relative standard deviation of the binomial frequency estimator is
	controlled by the probability of the selected event. A target relative standard
	deviation $\delta$ requires the indicative shot count~\cite{Nielsen2010}}
\begin{equation}
	\hll{S_{\mathrm{dir}}\simeq
		\frac{1-\Pr(\mathrm{success})}
		{\delta^2\Pr(\mathrm{success})}.}
	\label{eq:direct-postselection-shots}
\end{equation}
\hl{In the kinetic-energy construction of
	Section~\ref{sec:homogenized-observable}, the success probability of the solution branch is the already defined $P_K$, see Eq.~\eqref{eq:kinetic-branch-probability}.}

\hl{Equation~\eqref{eq:generic-polynomial-encoding} implements the scalar
	velocity map through a small-angle rotation. Writing the normalized Green
	multiplier as $\Gamma(\bm k)$, the linear target used there is
	$\theta_{\rm lin}(\bm k)=\varepsilon_\Gamma\Gamma(\bm k)$,
	whose selected amplitude is interpreted through $\sin(\theta_{\rm lin})
	\simeq\theta_{\rm lin}$. The reported value of $\varepsilon_\Gamma$ controls this linearization, as stated with the scalar-block scales in
	Section~\ref{subsec:operator-complexity}. It also sets the probability scale of $P_K$ and can make direct sampling costly.
	This direct-sampling estimate provides the baseline repetition count for the two mitigation strategies considered below.
}
\subsection{Arcsine-angle encoding for increasing selected amplitude}

\hl{We propose an alternative encoding that uses the available target-qubit amplitude more directly: it approximates}
\begin{equation}
	\hll{\theta_{\arcsin}(\bm k)=
		\arcsin\!\left(\rho_\Gamma\Gamma(\bm k)\right).}
	\label{eq:arcsin-gamma-encoding}
\end{equation}
\hl{We choose $\rho_\Gamma=0.8$ to use a large selected amplitude while
	keeping the argument of $\arcsin$ strictly inside its admissible interval and
	away from the rapidly increasing slope near its endpoints. The circuit applies
	$R_Y(2\theta)$ to the target qubit after approximating the corresponding
	linear or arcsine target angle with the piecewise-polynomial construction
	of Eq.~\eqref{eq:tile-polynomial}. 
	Assuming that the two angular fits have comparable error,
	the selected-branch probability scales can be compared as}
\begin{equation}
\hll{P_{K,{\rm lin}}\simeq\frac{2\varepsilon_\Gamma^2K}{F^2},
	\qquad
	P_{K,{\arcsin}}=
	\frac{2\rho_\Gamma^2K}{F^2}.}
\label{eq:linear-arcsin-probability-scales}
\end{equation}
\hl{At a fixed resolution, the two target functions correspond to the same
discrete kinetic energy $K$ and forcing norm $F^2$. Their probability ratio is
therefore}
\begin{equation}
\hll{\frac{P_{K,{\arcsin}}}{P_{K,{\rm lin}}}
	=\left(\frac{\rho_\Gamma}{\varepsilon_\Gamma}\right)^2.}
\label{eq:arcsin-probability-gain}
\end{equation}
\hl{With $\varepsilon_\Gamma=10^{-2}$ and $\rho_\Gamma=0.8$, the gain is
$(\rho_\Gamma/\varepsilon_\Gamma)^2=6400$. Polynomial approximation changes the recovered values slightly, as the RVE diagnostic verifies.}

\begin{figure}[tb]
	\centering
	\includegraphics[width=0.9\textwidth]{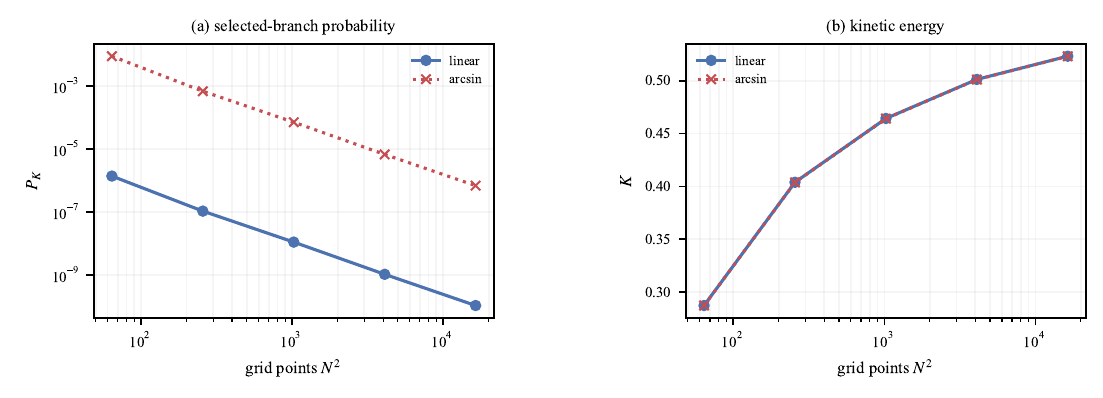}
	\caption{\label{fig:postselection-rve-mitigation}\hl{ RVE comparison of linear and arcsin encodings for the $U_\Gamma$ unitary. Panel (a) shows the selected-branch probability $P_K$; panel (b) shows the corresponding kinetic-energy predictions.}}
\end{figure}

\hl{Figure~\ref{fig:postselection-rve-mitigation} uses the RVE forcing of
Section~\ref{sec:results-homogenized}, with $U_R$ held at $p_R=3$.
Both encodings use the same degree-three monomial approximation on the same
eight-axis nested low-wavenumber tiling; only the fitted target changes from
$\varepsilon_\Gamma\Gamma$ to $\arcsin[\rho_\Gamma\Gamma]$. The dotted curves
reproduce the linear-angle encoding used in the main results, whereas the
solid curves use Eq.~\eqref{eq:arcsin-gamma-encoding}. The latter raises $P_K$
by approximately $6400$, while panel~(b) shows the corresponding kinetic-energy
predictions. The dipole calculation below clarifies the apparent decrease of
$P_K$ with $N$ for this RVE forcing.}

\begin{figure}[tb]
	\centering
	\includegraphics[width=0.9\textwidth]{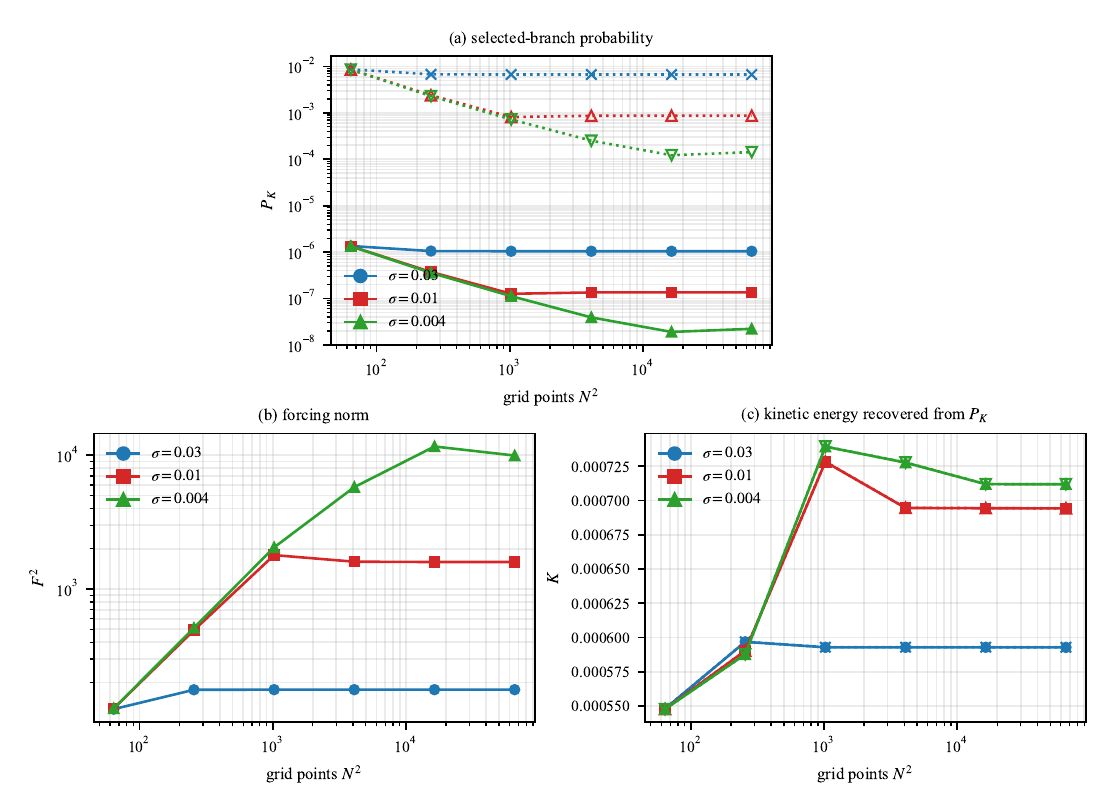}
	\caption{\label{fig:postselection-dipole-diagnostic}\hl{ Postselection
			refinement diagnostic for the regularized dipole. In panel (a), solid curves with filled markers
			use the linear small-angle encoding, while dotted curves use the arcsine-angle
			encoding with $\rho_\Gamma=0.8$. Panels (b) and (c) show the corresponding
			discrete forcing norm and the kinetic energy recovered from $P_K$ at the
			same grid resolution, respectively.}}
\end{figure}

\hl{In Figure~\ref{fig:postselection-dipole-diagnostic}, the regularized dipole shows that a decrease of $P_K$ under refinement is
not determined by the number of grid labels alone. Smaller $\sigma$ denotes a
narrower, more localized force and hence finer physical features. As refinement
resolves additional spectral content, $F^2$ can increase while the inverse
Stokes operator strongly attenuates the corresponding high-frequency
contribution to the kinetic-energy observable recovered from the selected
probability. Thus, $K$ can become nearly resolution-independent while $F^2$
continues to increase, and $P_K\propto K/F^2$ correspondingly decreases until
the relevant forcing content is resolved. Before this regime, the collocated
dipole samples are not nested spectral truncations, so aliasing can make $K$
non-monotone; once the active spectrum is resolved, it converges. The
arcsine-angle curves retain this resolution trend while increasing the overall
probability scale.}

\subsection{Amplitude amplification with observable recovery}
		
		\hl{Amplitude amplification coherently transfers amplitude from the
unselected branches to a chosen measurement event, increasing its probability
while retaining the mode registers coherently~\cite{Brassard2002,Grinko2021}. In the
present construction, the selected event is already specified by $\Pi_K$ in
Eq.~\eqref{eq:kinetic-branch-probability}. The two reflections required for
one amplification iteration are a phase flip of this selected event and a
phase flip of the input state:}
\begin{equation}
\hll{S_{\mathrm{good}}=I-2\Pi_K,\qquad
	S_0=I-2|0\rangle\langle0|,\qquad
	Q=-\mathcal A_K S_0\mathcal A_K^\dagger S_{\mathrm{good}}.}
\label{eq:aa-iterate}
\end{equation}
\hl{The reflection $S_{\mathrm{good}}$ acts only on the component, target,
and work registers, while the mode registers remain coherent. Writing
$P_K=q_0=\sin^2\vartheta$ for the probability before amplification, let $q_j$
denote the probability of the same selected event after $j$ applications of
$Q$. Thus}
\begin{equation}
\hll{q_j=\sin^2\!\left[(2j+1)\vartheta\right].}
\label{eq:aa-probability}
\end{equation}
\hl{For absolute-observable estimation, $j$ is chosen so that $q_j$ is close
to one half and $(2j+1)\vartheta\leq\pi/2$, so that the map from $P_K$ to
$q_j$ remains one-to-one. A measurement of $q_j$ then recovers the initial
probability $P_K=q_0$ through}
\begin{equation}
\hll{P_K=\sin^2\!\left[
	\frac{\arcsin\sqrt{q_j}}{2j+1}
	\right].}
\label{eq:aa-probability-inverse}
\end{equation}
\hl{The original kinetic energy then follows from
Eq.~\eqref{eq:kinetic-energy-from-probability}. For $S_{\rm AA}$ independent
executions of the amplified circuit, the sampled frequency of the selected
event estimates $q_j$ with standard deviation
$\sqrt{q_j(1-q_j)/S_{\rm AA}}$. Requiring a relative standard deviation
$\delta$ for this recovered probability gives}
\begin{equation}
\hll{S_{\rm AA}\simeq
	\frac{1-P_K}{M^2\delta^2P_K}.}
\label{eq:aa-postselection-shots}
\end{equation}
\hl{Table~\ref{tab:postselection-dipole-aa} applies this estimate to the
arcsine-angle dipole results at $N=256$. Here $S_{\rm dir}$ is the number of
direct circuit executions and $S_{\rm AA}$ is the number of executions of the amplified circuit. Let $j$ denote the iteration count and $M=2j+1$ the amplification multiplier; the column
$M S_{\rm AA}$ reports the resulting number of calls to the observable block
across all amplified executions.}

\begin{table}[tb!]
\centering
\caption{\label{tab:postselection-dipole-aa}\hl{ Ideal amplitude-amplification
		shot estimates for the arcsine-angle dipole calculation at $N=256$,
		$\rho_\Gamma=0.8$, and $\delta=10\%$. At $\delta=1\%$, every shot count is
		larger by a factor of $100$.}}
\begin{tabular}{crrrrrr}
	\toprule
	\(\sigma\) & \(P_K\) & \(M\) & \(q_j\) & \(S_{\rm dir}\) & \(S_{\rm AA}\) & \(M S_{\rm AA}\)\\
	\midrule
	\hll{0.03}  & \hll{6.69\!\times\!10^{-3}} & \hll{9}  & \hll{0.4516} & \hll{14\,849}  & \hll{184} & \hll{1\,656}\\
	\hll{0.01}  & \hll{8.71\!\times\!10^{-4}} & \hll{27} & \hll{0.5114} & \hll{114\,766} & \hll{158} & \hll{4\,266}\\
	\hll{0.004} & \hll{1.43\!\times\!10^{-4}} & \hll{65} & \hll{0.4914} & \hll{700\,027} & \hll{166} & \hll{10\,790}\\
	\bottomrule
\end{tabular}
\end{table}
\hl{For the most localized dipole, $\sigma=0.004$ at $N=256$, direct sampling
with the linear-angle encoding requires approximately $4.5\times10^9$ shots for
$\delta=10\%$. The arcsine-angle encoding reduces this to $7.0\times10^5$
direct shots, and ideal amplitude amplification with $M=65$ further gives
about $166$ amplified-circuit executions, corresponding to $1.08\times10^4$
observable-block calls. The reduction from the arcsine-angle direct-sampling
count to the amplified-circuit count is exchanged for coherent circuit depth.
One amplified execution contains $j+1$ applications of $\mathcal A_K$, $j$
applications of $\mathcal A_K^\dagger$, and $j$ implementations of each
reflection $S_{\rm good}$ and $S_0$; it therefore contains $M$ calls to the
observable block. These additional sequential operations require the quantum
state to remain coherent throughout the amplified-circuit execution. A complete resource assessment must combine
reflection synthesis, the procedure used to select $j$, state preparation,
hardware noise, fault tolerance, and the cost of broader multiscale readout.}

\bibliographystyle{unsrt}
\bibliography{biblio2}

@article{Succi2023,
  title     = {Quantum Computing for Fluids: Where Do We Stand?},
  author    = {Succi, Sauro and Itani, Wael and Sreenivasan, Katepalli and Steijl, Ren{\'e}},
  journal   = {Europhysics Letters},
  volume    = {144},
  number    = {1},
  pages     = {10001},
  year      = {2023},
  doi       = {10.1209/0295-5075/acfd19}
 }

@article{Turro2025,
  author  = {Turro, Francesco and Lignarolo, Alessandra and Dragoni, Daniele},
  title   = {Toward Practical Application of the Quantum {Carleman} {Lattice Boltzmann} Method in Industrial {CFD} Simulations},
  journal = {IEEE Transactions on Quantum Engineering},
  volume  = {6},
  pages   = {1--20},
  year    = {2025},
  doi     = {10.1109/TQE.2025.3620130},
  eprint  = {2504.13033},
  archivePrefix = {arXiv}
 }

@article{Febrianto2026,
  author  = {Febrianto, Eky and Wang, Yiren and Liu, Burigede and Ortiz, Michael and Cirak, Fehmi},
  title   = {A Quantum Spectral Method for Non-Periodic Boundary Value Problems},
  journal = {Computer Methods in Applied Mechanics and Engineering},
  volume  = {457},
  pages   = {118934},
  year    = {2026},
  doi     = {10.1016/j.cma.2026.118934},
  eprint  = {2511.11494},
  archivePrefix = {arXiv}
 }

@article{Raisuddin2026,
  author  = {Raisuddin, Osama Muhammad and De, Suvranu},
  title   = {A Review of Quantum Scientific Computing Algorithms Relevant to Computational Mechanics},
  journal = {Archives of Computational Methods in Engineering},
  volume  = {33},
  number  = {1},
  pages   = {745--787},
  year    = {2026},
  doi     = {10.1007/s11831-025-10321-9}
 }

@article{Pfeffer2025,
  author  = {Pfeffer, Philipp and Brearley, Peter and Laizet, Sylvain and Schumacher, J{\"o}rg},
  title   = {Spectral Quantum Algorithm for Passive Scalar Transport in Shear Flows},
  journal = {Scientific Reports},
  volume  = {15},
  number  = {41172},
  year    = {2025},
  doi     = {10.1038/s41598-025-27219-y},
  eprint  = {2505.10136},
  archivePrefix = {arXiv}
 }

@article{LiuStokes2024,
  author  = {Liu, Yunya and Chen, Zhen and Shu, Chang and Rebentrost, Patrick and Liu, Yangang and Chew, Shi Chuan and Cui, Yongdong},
  title   = {A Variational Quantum Algorithm-Based Numerical Method for Solving Potential and {Stokes} Flows},
  journal = {Ocean Engineering},
  volume  = {292},
  pages   = {116494},
  year    = {2024},
  doi     = {10.1016/j.oceaneng.2023.116494}
 }

@article{Liu2024,
  author  = {Liu, Burigede and Ortiz, Michael and Cirak, Fehmi},
  title   = {Towards Quantum Computational Mechanics},
  journal = {Computer Methods in Applied Mechanics and Engineering},
  volume  = {432},
  pages   = {117403},
  year    = {2024},
  doi     = {10.1016/j.cma.2024.117403},
  eprint  = {2312.03791},
  archivePrefix = {arXiv}
 }

@article{Idelsohn2024PDNS,
  author  = {Idelsohn, S. R. and Gimenez, J. M. and Larreteguy, A. E. and Nigro, N. M. and S{\'i}vori, F. M. and O{\~n}ate, E.},
  title   = {The {P-DNS} Method for Turbulent Fluid Flows: An Overview},
  journal = {Archives of Computational Methods in Engineering},
  volume  = {31},
  number  = {2},
  pages   = {973--1021},
  year    = {2024},
  doi     = {10.1007/s11831-023-09959-6}
 }

@article{Gimenez2025PDNS,
  author  = {Gimenez, J. M. and Sivori, F. M. and Larreteguy, A. E. and Montano, S. I. and Aguerre, H. J. and Orbaiz, P. and Idelsohn, S. R.},
  title   = {A Multiscale {Pseudo-DNS} Approach for Solving Turbulent Boundary-Layer Problems},
  journal = {Computer Methods in Applied Mechanics and Engineering},
  volume  = {437},
  pages   = {117804},
  year    = {2025},
  doi     = {10.1016/j.cma.2025.117804}
 }

@article{Wang2026,
      title={{QAFE$^2$}: Quantum accelerated multiscale finite element analysis},
      author={Wang, Yiren and Ortiz, Michael and Cirak, Fehmi},
      journal={Computer Methods in Applied Mechanics and Engineering},
      volume={461},
      pages={119204},
      year={2026},
      publisher={Elsevier}
    }

@article{GonzalezConde2025,
  author  = {Gonzalez-Conde, Javier and Lewis, Dylan and Bharadwaj, Sachin S. and Sanz, Mikel},
  title   = {Quantum State Preparation for Multivariate Functions},
  journal = {Quantum},
  volume  = {9},
  pages   = {1703},
  year    = {2025},
  doi     = {10.22331/q-2025-04-09-1703}
 }

@article{song25,
  author  = {Song, Zhixin and Deaton, Robert and Gard, Bryan T. and Bryngelson, Spencer H.},
  title   = {Incompressible {Navier--Stokes} Solve on Noisy Quantum Hardware via a Hybrid Quantum--Classical Scheme},
  journal = {Computers \& Fluids},
  volume  = {288},
  pages   = {106507},
  year    = {2025},
  doi     = {10.1016/j.compfluid.2024.106507},
  eprint  = {2406.00280},
  archivePrefix = {arXiv}
 }

@article{Joven2025,
  author  = {Joven, Kevin J. and Das, Elin Ranjan and Bierman, Joel and Majumdar, Aishwarya and Heris, Masoud Hakimi and Liu, Yuan},
  title   = {Scalable Quantum Computational Science: A Perspective from Block-Encodings and Polynomial Transformations},
  journal = {APL Computational Physics},
  volume  = {2},
  number  = {1},
  pages   = {010901},
  year    = {2026},
  doi     = {10.1063/5.0312254},
  eprint  = {2511.16738},
  archivePrefix = {arXiv}
 }

@book{Nielsen2010,
  author    = {Nielsen, Michael A. and Chuang, Isaac L.},
  title     = {Quantum Computation and Quantum Information},
  publisher = {Cambridge University Press},
  edition   = {10th anniversary},
  year      = {2010}
 }

@article{PleschBrukner2011,
  author  = {Plesch, Martin and Brukner, {\v{C}}aslav},
  title   = {Quantum-state Preparation with Universal Gate Decompositions},
  journal = {Physical Review A},
  volume  = {83},
  number  = {3},
  pages   = {032302},
  year    = {2011},
  doi     = {10.1103/PhysRevA.83.032302}
 }

@article{Sanchez2023,
  author  = {Marin-Sanchez, Guillermo and Gonzalez-Conde, Javier and Sanz, Mikel},
  title   = {Quantum Algorithms for Approximate Function Loading},
  journal = {Physical Review Research},
  volume  = {5},
  number  = {3},
  pages   = {033114},
  year    = {2023},
  doi     = {10.1103/PhysRevResearch.5.033114}
 }

@article{Zhang2022,
  author  = {Zhang, Xiao-Ming and Li, Tongyang and Yuan, Xiao},
  title   = {Quantum State Preparation with Optimal Circuit Depth: Implementations and Applications},
  journal = {Physical Review Letters},
  volume  = {129},
  number  = {23},
  pages   = {230504},
  year    = {2022},
  doi     = {10.1103/PhysRevLett.129.230504}
 }

@article{Qiskit2024,
  author        = {Javadi-Abhari, Ali and Treinish, Matthew and Krsulich, Kevin and Wood, Christopher J. and Lishman, Jake and Gacon, Julien and Martiel, Simon and Nation, Paul D. and Bishop, Lev S. and Cross, Andrew W. and Johnson, Blake R. and Gambetta, Jay M.},
  title         = {Quantum Computing with {Qiskit}},
  journal       = {arXiv preprint},
  eprint        = {2405.08810},
  archivePrefix = {arXiv},
  year          = {2024},
  doi           = {10.48550/arXiv.2405.08810}
 }

@article{schalkers2024momentum,
        title={Momentum exchange method for quantum Boltzmann methods},
        author={Schalkers, Merel A and M{\"o}ller, Matthias},
        journal={Computers \& Fluids},
        volume={285},
        pages={106453},
        year={2024},
        publisher={Elsevier}
    }

@article{georgescu2026efficient,
      title={Efficient and Expressive Boundary Conditions in Quantum Lattice Boltzmann Methods},
      author={Georgescu, C{\u{a}}lin A and M{\"o}ller, Matthias},
      journal={arXiv preprint arXiv:2606.01426},
      year={2026}
    }

@article{luacuatucs2026surrogate,
      title={Surrogate Quantum Circuit Design for the Lattice Boltzmann Collision Operator},
      author={L{\u{a}}c{\u{a}}tu{\c{s}}, Monica and M{\"o}ller, Matthias},
      journal={International Journal for Numerical Methods in Engineering},
      volume={127},
      number={4},
      pages={e70286},
      year={2026},
      publisher={Wiley Online Library}
    }

@article{tennie2025quantum,
      title={Quantum computing for nonlinear differential equations and turbulence},
      author={Tennie, Felix and Laizet, Sylvain and Lloyd, Seth and Magri, Luca},
      journal={Nature Reviews Physics},
      volume={7},
      number={4},
      pages={220--230},
      year={2025},
      publisher={Nature Publishing Group UK London}
    }

@article{Chen2024,
     title={Enabling large-scale and high-precision fluid simulations on near-term quantum computers},
     author={Chen, Zhao-Yun and Ma, Teng-Yang and Ye, Chuang-Chao and Xu, Liang and Bai, Wen and Zhou, Lei and Tan, Ming-Yang and Zhuang, Xi-Ning and Xu, Xiao-Fan and Wang, Yun-Jie and others},
     journal={Computer Methods in Applied Mechanics and Engineering},
     volume={432},
     pages={117428},
     year={2024},
     publisher={Elsevier}
    }

@article{Xu2024,
     title = {Quantum computing enhanced distance-minimizing data-driven computational mechanics},
     journal = {Computer Methods in Applied Mechanics and Engineering},
     volume = {419},
     pages = {116675},
     year = {2024},
     issn = {0045-7825},
     author = {Yongchun Xu and Jie Yang and Zengtao Kuang and Qun Huang and Wei Huang and Heng Hu}
    }

@article{Xu2025,
     title={Decomposition-free variational quantum linear solver: Application in computational mechanics},
     author={Xu, Yongchun and Hu, Heng},
     journal={Computer Methods in Applied Mechanics and Engineering},
     volume={447},
     pages={118396},
     year={2025},
     publisher={Elsevier}
    }

@article{Siegl2026,
     title={Tensor-programmable quantum circuits for solving differential equations},
     author={Siegl, Pia and Reese, Greta Sophie and Hashizume, Tomohiro and van H{\"u}lst, Nis-Luca and Jaksch, Dieter},
     journal={Physical Review Research},
     volume={8},
     number={1},
     pages={013052},
     year={2026},
     publisher={APS}
    }

@article{Novikau2025,
     title={Quantum algorithm for the advection-diffusion equation and the Koopman-von Neumann approach to nonlinear dynamical systems},
     author={Novikau, Ivan and Joseph, Ilon},
     journal={Computer Physics Communications},
     volume={309},
     pages={109498},
     year={2025},
     publisher={Elsevier}
    }

@article{Xue2025,
     title={Quantum homotopy analysis method with quantum-compatible linearization for nonlinear partial differential equations},
     author={Xue, Cheng and Xu, Xiao-Fan and Zhuang, Xi-Ning and Sun, Tai-Ping and Wang, Yun-Jie and Tan, Ming-Yang and Ye, Chuang-Chao and Liu, Huan-Yu and Wu, Yu-Chun and Chen, Zhao-Yun and others},
     journal={Science China Physics, Mechanics \& Astronomy},
     volume={68},
     number={10},
     pages={104702},
     year={2025},
     publisher={Springer}
    }

@book{Canuto1988,
      author    = {Canuto, Claudio and Hussaini, M. Yousuff and Quarteroni, Alfio and Zang, Thomas A.},
      title     = {Spectral Methods in Fluid Dynamics},
      publisher = {Springer},
      year      = {1988}
    }

@book{Trefethen2000,
      author    = {Trefethen, Lloyd N.},
      title     = {Spectral Methods in {MATLAB}},
      publisher = {SIAM},
      year      = {2000}
    }

@article{HHL2009,
      author  = {Harrow, Aram W. and Hassidim, Avinatan and Lloyd, Seth},
      title   = {Quantum Algorithm for Linear Systems of Equations},
      journal = {Physical Review Letters},
      volume  = {103},
      pages   = {150502},
      year    = {2009},
      doi     = {10.1103/PhysRevLett.103.150502}
    }

@inproceedings{Gilyen2019,
      author    = {Gily{\'e}n, Andr{\'a}s and Su, Yuan and Low, Guang Hao and Wiebe, Nathan},
      title     = {Quantum Singular Value Transformation and Beyond: Exponential Improvements for Quantum Matrix Arithmetics},
      booktitle = {Proceedings of the 51st Annual ACM SIGACT Symposium on Theory of Computing},
      pages     = {193--204},
      year      = {2019},
      doi       = {10.1145/3313276.3316366}
    }

@article{LowChuang2019,
      author  = {Low, Guang Hao and Chuang, Isaac L.},
      title   = {Hamiltonian Simulation by Qubitization},
      journal = {Quantum},
      volume  = {3},
      pages   = {163},
      year    = {2019},
      doi     = {10.22331/q-2019-07-12-163}
    }

@incollection{Brassard2002,
      author    = {Brassard, Gilles and H{\o}yer, Peter and Mosca, Michele and Tapp, Alain},
      title     = {Quantum Amplitude Amplification and Estimation},
      booktitle = {Quantum Computation and Information},
      series    = {Contemporary Mathematics},
      volume    = {305},
      pages     = {53--74},
      publisher = {American Mathematical Society},
      year      = {2002},
      doi       = {10.1090/conm/305/05215}
    }

@article{Grinko2021,
      author  = {Grinko, Dmitry and Gacon, Julien and Zoufal, Christa and Woerner, Stefan},
      title   = {Iterative Quantum Amplitude Estimation},
      journal = {npj Quantum Information},
      volume  = {7},
      pages   = {52},
      year    = {2021},
      doi     = {10.1038/s41534-021-00379-1}
    }

\end{document}